
%
%
%
%
%
%
\magnification=\magstephalf      
%
%
\vsize=7.5truein                 
\hsize=5.2truein                 
\newskip\stdskip                 
\stdskip=6pt plus3pt minus3pt    
\medskipamount=\stdskip          
\parindent=0pt                   
\parskip=\stdskip                
\abovedisplayskip=\stdskip       
\belowdisplayskip=\stdskip       
\mathsurround=0.75pt             
\overfullrule=0pt                
%
%
\def\ppar{\par\goodbreak\vskip 8pt plus 4pt minus 4pt}     
%
%
\def\stdspace{\hskip 0.75em plus 0.15em\ignorespaces}
\let\qua\stdspace 
%
%
%
%
%
%
%
\def\hexnumber#1{\ifcase#1 0\or 1\or 2\or 3\or 4\or 5\or 6\or 7\or 8\or
 9\or A\or B\or C\or D\or E\or F\fi}
%
%
\font\thirtnmsa=msam10 scaled 1315    
\font\tenmsa=msam10          \font\ninemsa=msam9
\font\sevenmsa=msam7         \font\sixmsa=msam6
\font\fivemsa=msam5
%
%
\newfam\msafam                  \textfont\msafam=\tenmsa
\scriptfont\msafam=\sevenmsa    \scriptscriptfont\msafam=\fivemsa
\edef\hexa{\hexnumber\msafam}        
\def\msa{\fam\msafam\tenmsa}         
%
%
\font\thirtnmsb=msbm10 scaled 1315   
\font\tenmsb=msbm10      \font\ninemsb=msbm9
\font\sevenmsb=msbm7     \font\sixmsb=msbm6
\font\fivemsb=msbm5
%
\newfam\msbfam                   \textfont\msbfam=\tenmsb       
\scriptfont\msbfam=\sevenmsb     \scriptscriptfont\msbfam=\fivemsb
\edef\hexb{\hexnumber\msbfam}    
\def\msb{\fam\msbfam\tenmsb}     
%
%
\font\thirtneufm=eufm10 scaled 1315   
\font\teneufm=eufm10                 \font\nineeufm=eufm9
\font\seveneufm=eufm7                \font\sixeufm=eufm6
\font\fiveeufm=eufm5
%
\newfam\eufmfam                    \textfont\eufmfam=\teneufm
\scriptfont\eufmfam=\seveneufm     \scriptscriptfont\eufmfam=\fiveeufm
\edef\hexf{\hexnumber\eufmfam}      
\def\frak{\fam\eufmfam\teneufm}     
%
%
%
\font\thirtnrm=cmr10 scaled 1315    
\font\ninerm=cmr9                   \font\sixrm=cmr6   
%
\font\thirtni=cmmi10 scaled 1315    
\font\ninei=cmmi9                   \font\sixi=cmmi6  
%
\font\thirtnsy=cmsy10 scaled 1315   
\font\ninesy=cmsy9                  \font\sixsy=cmsy6  
%
\font\thirtnbf=cmbx10 scaled 1315   
\font\ninebf=cmbx9                  \font\sixbf=cmbx6  
%
%
\font\thirtnex=cmex10 scaled 1315   
\font\nineex=cmex9                  
%
%
\font\thirtnit=cmti10 scaled 1315  
\font\nineit=cmti9                  
%
\font\thirtnsl=cmsl10 scaled 1315  
\font\ninesl=cmsl9                  
%
\font\thirtntt=cmtt10 scaled 1315  
\font\ninett=cmtt9                  
%
%
%
%
\def\small{%
%
%
\textfont0=\ninerm \scriptfont0=\sixrm \scriptscriptfont0=\fiverm
\def\rm{\fam0\ninerm}
%
%
\textfont1=\ninei \scriptfont1=\sixi \scriptscriptfont1=\fivei
%
%
\textfont2=\ninesy \scriptfont2=\sixsy \scriptscriptfont2=\fivesy
%
%
\textfont3=\nineex \scriptfont3=\nineex \scriptscriptfont3=\nineex
%
%
\textfont\bffam=\ninebf \scriptfont\bffam=\sixbf
\scriptscriptfont\bffam=\fivebf \def\bf{\fam\bffam\ninebf}%
%
%
\textfont\itfam=\nineit \def\it{\fam\itfam\nineit}%
\textfont\slfam=\ninesl \def\sl{\fam\slfam\ninesl}%
\textfont\ttfam=\ninett \def\tt{\fam\ttfam\ninett}%
%
%
%
\textfont\msafam=\ninemsa \scriptfont\msafam=\sixmsa
\scriptscriptfont\msafam=\fivemsa \def\msa{\fam\msafam\ninemsa}%
%
%
\textfont\msbfam=\ninemsb \scriptfont\msbfam=\sixmsb
\scriptscriptfont\msbfam=\fivemsb \def\msb{\fam\msbfam\ninemsb}%
%
%
\textfont\eufmfam=\nineeufm  \scriptfont\eufmfam=\sixeufm
\scriptscriptfont\eufmfam=\fiveeufm \def\frak{\fam\eufmfam\nineeufm}%
%
%
%
\normalbaselineskip=11pt%
\setbox\strutbox=\hbox{\vrule height8pt depth3pt width0pt}%
%
%
\normalbaselines\rm
%
%
\stdskip=4pt plus2pt minus2pt    
\medskipamount=\stdskip          
\parskip=\stdskip                
\abovedisplayskip=\stdskip       
\belowdisplayskip=\stdskip       
\def\ppar{\par\goodbreak\vskip 6pt plus 3pt minus 3pt}%
%
%
\def\section##1{\global\advance\sectionnumber by 1
\vskip-\lastskip\penalty-800\vskip 20pt plus10pt minus5pt 
\egroup{\bf\number\sectionnumber\quad##1}\bgroup\small         
\vskip 6pt plus3pt minus3pt
\nobreak\resultnumber=1}
}    
%
\def\beginsmall{\bgroup\small}
\let\endsmall\egroup
%
%
%
%
\def\large{%
\textfont0=\thirtnrm \scriptfont0=\ninerm \scriptscriptfont0=\sevenrm
\def\rm{\fam0\thirtnrm}%
\textfont1=\thirtni \scriptfont1=\ninei \scriptscriptfont1=\seveni
\textfont2=\thirtnsy \scriptfont2=\ninesy \scriptscriptfont2=\sevensy
\textfont3=\thirtnex \scriptfont3=\thirtnex \scriptscriptfont3=\thirtnex
\textfont\bffam=\thirtnbf \scriptfont\bffam=\ninebf
\scriptscriptfont\bffam=\sevenbf \def\bf{\fam\bffam\thirtnbf}%
\textfont\itfam=\thirtnit \def\it{\fam\itfam\thirtnit}%
\textfont\slfam=\thirtnsl \def\sl{\fam\slfam\thirtnsl}%
\textfont\ttfam=\thirtntt \def\tt{\fam\ttfam\thirtntt}%
\textfont\msafam=\thirtnmsa \scriptfont\msafam=\ninemsa
\scriptscriptfont\msafam=\sevenmsa \def\msa{\fam\msafam\thirtnmsa}%
\textfont\msbfam=\thirtnmsb \scriptfont\msbfam=\ninemsb
\scriptscriptfont\msbfam=\sevenmsb \def\msb{\fam\msbfam\thirtnmsb}%
\textfont\eufmfam=\thirtneufm  \scriptfont\eufmfam=\nineeufm
\scriptscriptfont\eufmfam=\seveneufm \def\frak{\fam\eufmfam\teneufm}%
\normalbaselineskip=16pt%
\setbox\strutbox=\hbox{\vrule height11.5pt depth4.5pt width0pt}%
\normalbaselines\rm}%
\let\Large\large   
%
\def\Bbb#1{{\msb#1}}

%

\def\re{\Bbb R}
%
\mathchardef\plussquare="0\hexa01
\mathchardef\nge="3\hexb0B
\mathchardef\maltesecross="0\hexa7A
\mathchardef\del="0\hexf01
%
%
%
%
\font\sc=cmcsc10
%
%
%
%
\def\sqr#1#2{{\vcenter{\vbox{\hrule  height.#2truept
	\hbox{\vrule width.#2truept height#1truept 
	\kern#1truept \vrule width.#2truept}
	\hrule height.#2truept}}}}
\def\sq{\sqr55}    
%
%
%
%
\newcount\sectionnumber            
\newcount\resultnumber             
\sectionnumber=0\resultnumber=1    
%
%
%
\def\section#1{\global\advance\sectionnumber by 1
\xdef\nextkey{\number\sectionnumber}
\vskip-\lastskip\penalty-800\vskip 20pt plus10pt minus5pt 
{\large\bf\number\sectionnumber\quad#1}         
\vskip 8pt plus4pt minus4pt
\nobreak\resultnumber=1}      
%
%
%
%
%
\def\sh#1{\vskip-\lastskip\ppar{\bf #1}\par\nobreak\medskip}         
%
%
%
%

%
\def\proc#1{\xdef\nextkey{\number\sectionnumber.\number\resultnumber}%
\vskip-\lastskip\ppar\bf%
\noindent#1\ \number\sectionnumber.\number\resultnumber
\stdspace\sl\global\advance\resultnumber by 1\ignorespaces}
 
%
%
\def\prf{\vskip-\lastskip\ppar\noindent{\bf Proof}%
\stdspace\rm}                            
\def\qed{\hfill$\sq$\par\goodbreak\rm}   
\def\endprf{\unskip\stdspace\hbox{}
\hfill$\sq$\par\medskip}                 
%
%
%
%
%
%
%
%
\def\proclaim#1{\vskip-\lastskip\ppar\bf%
\noindent#1\stdspace\sl\ignorespaces} 
\let\endproclaim\endproc
%
%
%
%
\def\rk#1{\vskip-\lastskip\ppar{\bf #1}\stdspace\ignorespaces}                

%
%
%
%
%
%
\def\label{\xdef\nextkey{\number\sectionnumber.\number\resultnumber}%
\number\sectionnumber.\number\resultnumber
\global\advance\resultnumber by 1}
%
%
%
%
%
%
%
%
%
%
%
%
%
%
%
%
\newcount\refnumber              
\refnumber=1                     
\long\def\reflist#1\endreflist{%
\long\def\thereflist{#1}{\def\refkey##1##2\par{\xdef##1{\number\refnumber}%
\global\advance\refnumber by 1}%
\def\key##1##2\par{\expandafter\xdef%
\csname##1\endcsname{\number\refnumber}%
\global\advance\refnumber by 1}#1\par}}
\long\def\references{%
\penalty-800\vskip-\lastskip\vskip 15pt plus10pt minus5pt 
{\large\bf References}\ppar 
{\leftskip=25pt\frenchspacing    
\small\parskip=3pt plus2pt       
\def\refkey##1##2\par{\noindent  
\llap{[##1]\stdspace}\ignorespaces##2\par}         
\def\key##1##2\par{\noindent  
\llap{[\ref{##1}]\stdspace}\ignorespaces##2\par}  
\def\,{\thinspace}\thereflist\par}}
%
%
%
\newcount\footnotenumber         
\footnotenumber=1                
\def\fnote#1{\xdef\nextkey{\number\footnotenumber}%
{\small\ifnum\footnotenumber>9\parindent=14pt%
\else\parindent=10pt\fi\footnote{$^{\number\footnotenumber}$}%
{\hglue-5pt#1}\global\advance\footnotenumber by 1}}
%
%
%
%
%
%
%
\newcount\figurenumber          
\figurenumber=1                 
\def\caption#1{\xdef\nextkey{\number\figurenumber}%
\cl{\small Figure \number\figurenumber: #1}%
\global\advance\figurenumber by 1}
\def\figurelabel{\xdef\nextkey{\number\figurenumber}%
\cl{\small Figure \number\figurenumber}%
\global\advance\figurenumber by 1}
\long\def\figure#1\endfigure{{\xdef\nextkey{\number\figurenumber}%
\let\captiontext\relax\def\caption##1{\xdef\captiontext{##1}}%
\midinsert\cl{\ignorespaces#1\unskip\unskip\unskip\unskip}\vglue6pt\cl{\small 
Figure \number\figurenumber\ifx\captiontext\relax\else: \captiontext
\fi}\endinsert\global\advance\figurenumber by 1}}
%
%
%
%
%
%
%
\def\nextkey{??}   
%
\def\key#1{\expandafter\xdef\csname #1\endcsname{\nextkey}}
\def\ref#1{\expandafter\ifx\csname #1\endcsname\relax
\immediate\write16{Reference {#1} undefined}??\else
\csname #1\endcsname\fi}
%
%
%
%
%
%
%
\newread\gtinfile
\newwrite\gtreffile
\def\useforwardrefs{
\openin\gtinfile\jobname.ref
\ifeof\gtinfile
\closein\gtinfile
\immediate\write16{No file \jobname.ref}
\else
\closein\gtinfile
\input \jobname.ref
\fi
\immediate\openout\gtreffile \jobname.ref
%
%
\def\key##1{{\def\\{\noexpand}%
\expandafter\xdef\csname ##1\endcsname{\nextkey}%
\immediate\write\gtreffile{\\\expandafter\\\def\\\csname ##1\\\endcsname%
{\nextkey}}}}
%
%
\long\def\reflist##1\endreflist{%
\long\def\thereflist{##1}{\def\refkey####1####2\par{\xdef####1{%
\number\refnumber}{\def\\{\noexpand}\immediate\write\gtreffile
{\\\def\\####1{\number\refnumber}}}\global\advance\refnumber by 1}%
\def\key####1####2\par{\expandafter\xdef%
\csname####1\endcsname{\number\refnumber}%
{\def\\{\noexpand}\immediate\write\gtreffile
{\\\expandafter\\\def\\\csname ####1\\\endcsname{\number\refnumber}}}
\global\advance\refnumber by 1}##1\par}}
\long\def\biblio##1\endbiblio{\reflist##1\endreflist\references}%
%
%
\def\numkey##1{{\def\\{\noexpand}%
\xdef##1{\number\sectionnumber.\number\resultnumber}
\immediate\write\gtreffile{\\\def\\##1%
{\number\sectionnumber.\number\resultnumber}}}}
\def\seckey##1{{\def\\{\noexpand}\xdef##1{\number\sectionnumber}
\immediate\write\gtreffile{\\\def\\##1{\number\sectionnumber}}}}
\def\figkey##1{\xdef##1{\number\figurenumber}%
{\def\\{\noexpand}\immediate\write\gtreffile%
{\\\def\\##1{\number\figurenumber}}}
\number\figurenumber\global\advance\figurenumber by 1}
}   
%
%
%
%
\def\figkey#1{\xdef#1{\number\figurenumber}%
\number\figurenumber\global\advance\figurenumber by 1}
\def\fig#1#2\endfig{%
\midinsert\cl{#2}\vglue6pt\cl{\small Figure #1}\endinsert}
\def\newfig{\number\figurenumber\global\advance\figurenumber by 1}
\def\numkey#1{\xdef#1{\number\sectionnumber.\number\resultnumber}}
\def\seckey#1{\xdef#1{\number\sectionnumber}}
%
%
%
%
%
%
%
%
%
\def\verb{\catcode`\"=\active}       
\def\brev{\catcode`\"=12}            
\brev                                
\verb                                
{\obeyspaces\gdef {\ }}              
{\catcode`\`=\active\gdef`{\relax\lq}}
\def"{%
\begingroup\baselineskip=12pt\def\par{\leavevmode\endgraf}%
\tt\obeylines\obeyspaces\parskip=0pt\parindent=0pt%
\catcode`\$=12\catcode`\&=12\catcode`\^=12\catcode`\#=12%
\catcode`\_=12\catcode`\~=12%
\catcode`\{=12\catcode`\}=12\catcode`\%=12\catcode`\\=12%
\catcode`\`=\active\let"\endgroup}
\brev      
%
%
%
%
%
%
\def\enditems{\par\leftskip = 0pt}         
\def\item#1{\par\leavevmode\llap{#1\stdspace}%
\ignorespaces}                             
%
%

%
%
\def\co{\colon\thinspace}    
\def\np{\vfil\eject}         
\def\nl{\hfil\break}         
\def\cl{\centerline}         
\def\agt{{\mathsurround=0pt\it$\cal A\mskip-.7mu$lgebraic \&\ 
$\cal G\mskip-2mu$eometric $\cal T\!\!$opology}}  
%
%
%

%
%
%
%
%
\def\title#1{\def\thetitle{#1}}

\def\author#1{\edef\previousauthors{\theauthors}
 \ifx\theauthors\relax\def\theauthors{#1}\else
 \def\theauthors{\previousauthors\par#1}\fi}

%
\def\address#1{\edef\previousaddresses{\theaddress}
 \ifx\theaddress\relax\def\theaddress{#1}\else
 \def\theaddress{\previousaddresses\par\vskip 2pt\par#1}\fi}
\def\secondaddress#1{\edef\previousaddresses{\theaddress}
 \ifx\theaddress\relax\def\theaddress{#1}\else
 \def\theaddress{\previousaddresses\par{\rm and}\par#1}\fi}   

\def\email#1{\edef\previousemails{\theemail}
 \ifx\theemail\relax\def\theemail{#1}\else
 \def\theemail{\previousemails\hskip 0.75em\relax#1}\fi}
\def\secondemail#1{\edef\previousemails{\theemail}
 \ifx\theemail\relax\def\theemail{#1}\else
 \def\theemail{\previousemails\hskip 0.75em{\rm and}\hskip 0.75em
 \relax#1}\fi}
\def\url#1{\edef\previousurls{\theurl}
 \ifx\theurl\relax\def\theurl{#1}\else
 \def\theurl{\previousurls\hskip 0.75em\relax#1}\fi}
\def\secondurl#1{\edef\previousurls{\theurl}
 \ifx\theurl\relax\def\theurl{#1}\else
 \def\theurl{\previousurls\hskip 0.75em{\rm and}\hskip 0.75em
 \relax#1}\fi}
\long\def\abstract#1\endabstract{\long\def\theabstract{#1}}
\def\primaryclass#1{\def\theprimaryclass{#1}}

\def\keywords#1{\def\thekeywords{#1}}
%
%
\let\\\par\let\thetitle\relax\let\theshorttitle\relax
\let\theauthors\relax\let\theshortauthors\relax
\let\theaddress\relax\let\theshortaddress\relax
\let\theemail\relax\let\theurl\relax
\let\theabstract\relax\let\theprimaryclass\relax
\let\thesecondaryclass\relax\let\thekeywords\relax
%
%
%
%
\long\def\maketitlepage{    

\vglue 0.2truein   

%
{\parskip=0pt\leftskip 0pt plus 1fil\def\\{\par\smallskip}{\large
\bf\thetitle}\par\medskip}   

\vglue 0.15truein 

%
{\parskip=0pt\leftskip 0pt plus 1fil\def\\{\par}{\sc\theauthors}
\par\medskip}%
 
\vglue 0.1truein 

%
{\small\parskip=0pt
{\leftskip 0pt plus 1fil\def\\{\par}{\sl\theaddress}\par}
\ifx\theemail\relax\else  
\vglue 5pt \def\\{\stdspace{\rm and}\stdspace} 
\cl{Email:\stdspace\tt\theemail}\fi
\ifx\theurl\relax\else    
\vglue 5pt \def\\{\stdspace{\rm and}\stdspace} 
\cl{URL:\stdspace\tt\theurl}\fi\par}

\vglue 7pt 

{\bf Abstract}

\vglue 5pt

\theabstract

\vglue 7pt 

{\bf AMS Classification numbers}\quad Primary:\quad \theprimaryclass\par

Secondary:\quad \thesecondaryclass

\vglue 5pt 

{\bf Keywords:}\quad \thekeywords

\np  

}    
%
%
\long\def\makeshorttitle{    


%
{\parskip=0pt\leftskip 0pt plus 1fil\def\\{\par\smallskip}{\large
\bf\thetitle}\par\medskip}   

\vglue 0.05truein 

%
{\parskip=0pt\leftskip 0pt plus 1fil\def\\{\par}{\sc\theauthors}
\par\medskip}%
 
\vglue 0.03truein 

%
{\small\parskip=0pt
{\leftskip 0pt plus 1fil\def\\{\par}{\sl\ifx\theshortaddress\relax
\theaddress\else\theshortaddress\fi}\par}
\ifx\theemail\relax\else  
\vglue 5pt \def\\{\stdspace{\rm and}\stdspace} 
\cl{Email:\stdspace\tt\theemail}\fi
\ifx\theurl\relax\else    
\vglue 5pt \def\\{\stdspace{\rm and}\stdspace} 
\cl{URL:\stdspace\tt\theurl}\fi\par}

\vglue 10pt 


{\small\leftskip 25pt\rightskip 25pt{\bf Abstract}\stdspace\theabstract

{\bf AMS Classification}\stdspace\theprimaryclass
\ifx\thesecondaryclass\relax\else; \thesecondaryclass\fi\par
{\bf Keywords}\stdspace \thekeywords\par}
\vglue 7pt
}    
\let\maketitle\makeshorttitle        
%
%

\def\volumenumber#1{\def\thevolumenumber{#1}}
\def\volumeyear#1{\def\thevolumeyear{#1}}
\def\pagenumbers#1#2{\def\startpage{#1}\def\finishpage{#2}}
\def\published#1{\def\publishdate{#1}}
\def\received#1{\def\receiveddate{#1}}

\let\reviseddate\relax
\volumenumber{X}
\volumeyear{20XX}
\pagenumbers{1}{XXX}
\published{XX Xxxember 20XX}

\long\def\makeagttitle{   
\agt\hfill      
\hbox to 60truept{\vbox to 0pt{\vglue -14truept{\bf [Logo here]}\vss}\hss}
\break
{\small Volume \thevolumenumber\ (\thevolumeyear)
\startpage--\finishpage\nl
Published: \publishdate}

\vglue .2truein

{\parskip=0pt\leftskip 0pt plus 1fil\def\\{\par\smallskip}{\large
\bf\thetitle}\par\medskip}   
\vglue 0.05truein 

%
{\parskip=0pt\leftskip 0pt plus 1fil\def\\{\par}{\sc\theauthors}
\par\medskip}%
 
\vglue 0.03truein 


{\small\leftskip 25truept\rightskip 25truept{\bf Abstract}\stdspace\theabstract

{\bf AMS Classification}\stdspace\theprimaryclass
\ifx\thesecondaryclass\relax\else; \thesecondaryclass\fi\par
{\bf Keywords}\stdspace \thekeywords\par}\vglue 7truept

}   


\def\Addresses{\bigskip
{\small \parskip 0pt \leftskip 0pt \rightskip 0pt plus 1fil \def\\{\par}
\sl\theaddress\par\medskip \rm Email:\stdspace\tt\theemail\par
\ifx\theurl\relax\else\smallskip \rm URL:\stdspace\tt\theurl\par\fi}}

\def\agtart{
\hoffset 14truemm
\voffset 31truemm
\font\phead=cmsl9 scaled 950
\font\pnum=cmbx10 scaled 913
\font\pfoot=cmsl9 scaled 950
\headline{\vbox to 0pt{\vskip -4.5mm\line{\small\phead\ifnum
\count0=\startpage ISSN numbers are printed here
\hfill {\pnum\folio}\else\ifodd\count0\def\\{ }%
\ifx\theshorttitle\relax\thetitle\else\theshorttitle\fi\hfill{\pnum\folio}
\else\def\\{ and }{\pnum\folio}\hfill\ifx\theshortauthors\relax\theauthors
\else\theshortauthors\fi\fi\fi}\vss}}
\footline{\vbox to 0pt{\vglue 0mm\line{\small\pfoot\ifnum\count0=\startpage
Copyright declaration is printed here\hfill\else
\agt, Volume \thevolumenumber\ (\thevolumeyear)\hfill\fi}\vss}}
\let\maketitle\makeagttitle\let\makeshorttitle\makeagttitle}


\def\ifplaintex{\expandafter\ifx\csname documentclass\endcsname\relax}

\def\gtp{{\mathsurround=0pt\it $\cal G\mskip-2mu$eometry \&\ 
$\cal T\!\!$opology $\cal P\!$ublications}}  

\def\Addressesr{\bigskip
{\small \parskip 0pt \leftskip 0pt \rightskip 0pt plus 1fil \def\\{\par}
\sl\theaddress\par
\medskip
\rm Email:\stdspace\tt\theemail\hfill\rm Received:\qua\receiveddate \par}}

\def\recd{{\small Received:\qua\receiveddate\ifx\reviseddate\relax
\else\qquad Revised:\qua\reviseddate\fi\par}} 


\def\lognumber#1{\def\thelognumber{#1}}
\def\volumenumber#1{\def\thevolumenumber{#1}}
\def\volumeyear#1{\def\thevolumeyear{#1}}
\def\papernumber#1{\def\thepapernumber{#1}}
\def\pagenumbers#1#2{\def\startpage{#1}\def\finishpage{#2}}
\def\published#1{\def\publishdate{#1}}

\def\received#1{\def\receiveddate{#1}}

\def\accepted#1{\def\accepteddate{#1}}


\let\\\par\let\thelognumber\relax\let\thevolumenumber\relax
\let\thepapernumber\relax\let\thevolumeyear\relax\let\startpage\relax
\let\finishpage\relax\let\publishdate\relax\let\receiveddate\relax
\let\reviseddate\relax\let\accepteddate\relax\let\theasciititle\relax
\let\theasciiauthors\relax
\let\theasciiabstract\relax

\let\theasciiemail\relax


\ifplaintex
\font\logobig=cmssbx10 scaled 3836
\font\logomed=cmssbx10 scaled 2557
\else
\font\logobig=cmssbx10 scaled 4200
\font\logomed=cmssbx10 scaled 2800
\fi

\long\def\makeagttitle{   
\count0=\startpage
\agt\hfill      
\hbox to 45truept{\vbox to 0pt{\vglue -13truept{\logomed A\kern -.37em{\logobig 
T}\kern -.38em G}\vss}\hss}
\break
{\small Volume \thevolumenumber\ (\thevolumeyear)
\startpage--\finishpage\nl
Published: \publishdate}

\vglue .25truein

{\parskip=0pt\leftskip 0pt plus
1fil\def\\{\par\smallskip}{\Large\bf\thetitle}\par\medskip} \vglue
0.05truein

%
{\parskip=0pt\leftskip 0pt plus 1fil\def\\{\par}{\sc\theauthors}
\par\medskip}%
 
\vglue 0.03truein 


{\small\leftskip 25truept\rightskip 25truept{\bf Abstract}\stdspace\theabstract

{\bf AMS Classification}\stdspace\theprimaryclass
\ifx\thesecondaryclass\relax\else; \thesecondaryclass\fi\par
{\bf Keywords}\stdspace \thekeywords\par}\vglue 7truept

}   

\ifplaintex
\hoffset 14truemm
\voffset 31truemm
\font\phead=cmsl9 scaled 950
\font\pnum=cmbx10 scaled 913
\font\pfoot=cmsl9 scaled 950
\headline{\vbox to 0pt{\vskip -4.5mm\line{\small\phead\ifnum
\count0=\startpage ISSN 1472-2739 (on-line) 1472-2747 (printed)
\hfill {\pnum\folio}\else\ifodd\count0\def\\{ }%
\ifx\theshorttitle\relax\thetitle\else\theshorttitle\fi\hfill{\pnum\folio}
\else\def\\{ and }{\pnum\folio}\hfill\ifx\theshortauthors\relax\theauthors
\else\theshortauthors\fi\fi\fi}\vss}}
\footline{\vbox to 0pt{\vglue 0mm\line{\small\pfoot\ifnum\count0=\startpage
\copyright\ \gtp\hfill\else
\agt, Volume \thevolumenumber\ (\thevolumeyear)\hfill\fi}\vss}}
\else
\headsep 23pt
\footskip 35pt
\hoffset -4truemm
\voffset 12.5truemm
\font\lhead=cmsl9 scaled 1050
\font\lnum=cmbx10 
\font\lfoot=cmsl9 scaled 1050
\makeatletter
\def\@oddhead{{\small\lhead\ifnum\count0=\startpage ISSN 1472-2739 
(on-line) 1472-2747 (printed)\hfill {\lnum\number\count0}\else\ifodd\count0
\def\\{ }\ifx\theshorttitle\relax \thetitle \else\theshorttitle\fi\hfill
{\lnum\number\count0}\else\def\\{ and }{\lnum\number\count0}
\hfill\ifx\theshortauthors\relax 
\theauthors\else\theshortauthors\fi\fi\fi}}\def\@evenhead{\@oddhead}
\def\@oddfoot{\small\lfoot\ifnum\count0=\startpage\copyright\ \gtp\hfill\else
\agt, Volume \thevolumenumber\ (\thevolumeyear)\hfill\fi}
\def\@evenfoot{\@oddfoot}
\makeatother
\fi
\let\maketitlepage\makeagttitle
\let\makeshorttitle\maketitlepage
\let\maketitle\maketitlepage


\newwrite\gtoutfile
\long\gdef\makeheadfile{  
{\def\\{, }\def\s{ }
\immediate\openout\gtoutfile head.xxx
\immediate\write\gtoutfile{To: math@arxiv.org}
\immediate\write\gtoutfile{Subject: put OR rep NNNNN:ppppp}
\immediate\write\gtoutfile{--text follows this line--}
\immediate\write\gtoutfile{Proxy-for: \ifx\theasciiauthors\relax
\theauthors\else\theasciiauthors\fi\s<\ifx\theasciiemail\relax\theemail\else\theasciiemail\fi>}
\immediate\write\gtoutfile{\noexpand\\}
\immediate\write\gtoutfile{Authors: \ifx\theasciiauthors\relax
\theauthors\else\theasciiauthors\fi}
{\def\\{ }\immediate\write\gtoutfile{Title: \ifx\theasciititle\relax
\thetitle\else\theasciititle\fi}}
\immediate\write\gtoutfile{Subj-class: GT or SG, GR etc}
\immediate\write\gtoutfile{MSC-class: \theprimaryclass\ifx\thesecondaryclass\relax\else, \thesecondaryclass\fi}
\immediate\write\gtoutfile{Journal-ref: Algebr. Geom. Topol. \thevolumenumber\s
(\thevolumeyear) \startpage-\finishpage}
\immediate\write\gtoutfile{Comments: Published by Algebraic and
Geometric Topology at}
\immediate\write\gtoutfile{\s\s\s  http://www.maths.warwick.ac.uk/agt/AGTVol\thevolumenumber/agt-\thevolumenumber-\thepapernumber.abs.html}
\immediate\write\gtoutfile{\noexpand\\}
\immediate\write\gtoutfile{}
\ifx\theasciiabstract\relax
\immediate\write\gtoutfile{\theabstract}\else
\immediate\write\gtoutfile{\theasciiabstract}\fi
\immediate\write\gtoutfile{}
\immediate\write\gtoutfile{\noexpand\\}
\immediate\write\gtoutfile{}
\immediate\closeout\gtoutfile}}  

\def\maketitlepage{\makeagttitle\makeheadfile}
\let\makeshorttitle\maketitlepage
\let\maketitle\maketitlepage

\lognumber{25}
\volumenumber{3}
\volumeyear{2003}
\papernumber{25}
\published{3 August 2003}
\pagenumbers{719}{775}
\received{8 November 2002}
\accepted{13 March 2003}

\input pictex

\def\bigstrut{\vrule width 0pt height 15pt depth 10pt}
\def\smallstrut{\vrule width 0pt height 8pt depth 5pt}
\def\medstrut{\vrule width 0pt height 10pt depth 5pt}


\def\ltreea#1#2#3#4#5#6#7{
\beginpicture
\small
\setcoordinatesystem units <0.5cm,0.6cm>
\linethickness=1pt
\setplotsymbol ({\tenrm .})
\setlinear
\plot 0 -.25  0 1  -1 2  0 3 /
\plot 0 1  2 3 /
\plot -1 2  -2 3 /
\put {$#7$} [Br] <-3pt, 0pt> at 0 -.25
\put {$#1$} [b] <0pt, 3pt> at -2 3
\put {$#2$} [b] <0pt, 3pt> at 0 3
\put {$#3$} [b] <0pt, 3pt> at 2 3
\put {$#5$} [tr] <-3pt, 0pt> at -.5 1.5
\put {$#6$} [b] <0pt, 5pt> at 0 1
\put {$#4$} [b] <0pt, 5pt> at -1 2
\endpicture}

\def\rtreea#1#2#3#4#5#6#7{
\beginpicture
\small
\setcoordinatesystem units <-0.5cm,0.6cm>
\linethickness=1pt
\setplotsymbol ({\tenrm .})
\setlinear
\plot 0 -.25  0 1  -1 2  0 3 /
\plot 0 1  2 3 /
\plot -1 2  -2 3 /
\put {$#7$} [Br] <-3pt, 0pt> at 0 -.25
\put {$#3$} [b] <0pt, 3pt> at -2 3
\put {$#2$} [b] <0pt, 3pt> at 0 3
\put {$#1$} [b] <0pt, 3pt> at 2 3
\put {$#5$} [tl] <3pt, 0pt> at -.5 1.5
\put {$#6$} [b] <0pt, 5pt> at 0 1
\put {$#4$} [b] <0pt, 5pt> at -1 2
\endpicture}

\def\ltreeb#1#2#3#4#5#6#7#8{
\beginpicture
\small
\setcoordinatesystem units <0.5cm,0.6cm>
\linethickness=1pt
\setplotsymbol ({\tenrm .})
\setlinear
\plot 0 -0.5  0 1  -3 4 /
\plot 0 1  3 4 /
\plot -1 2  1 4 /
\plot -2 3  -1 4 /
\put {$#1$} [b] <0pt, 3pt> at -3 4
\put {$#1$} [b] <0pt, 3pt> at -1 4
\put {$#1$} [b] <0pt, 3pt> at 1 4
\put {$#1$} [b] <0pt, 3pt> at 3 4
\put {$#2$} [b] <0pt, 5pt> at -2 3
\put {$#3$} [tr] <-3pt, 0pt> at -2.5 3.5
\put {$#4$} [b] <0pt, 5pt> at -1 2
\put {$#5$} [tr] <-3pt, 0pt> at -1.5 2.5
\put {$#6$} [b] <0pt, 5pt> at 0 1
\put {$#7$} [tr] <-3pt, 0pt> at -0.5 1.5
\put {$#8$} [Br] <-3pt, 0pt> at 0 -0.5
\endpicture}

\def\rtreeb#1#2#3#4#5#6#7#8{
\beginpicture
\small
\setcoordinatesystem units <-0.5cm,0.6cm>
\linethickness=1pt
\setplotsymbol ({\tenrm .})
\setlinear
\plot 0 -0.5  0 1  -3 4 /
\plot 0 1  3 4 /
\plot -1 2  1 4 /
\plot -2 3  -1 4 /
\put {$#1$} [b] <0pt, 3pt> at -3 4
\put {$#1$} [b] <0pt, 3pt> at -1 4
\put {$#1$} [b] <0pt, 3pt> at 1 4
\put {$#1$} [b] <0pt, 3pt> at 3 4
\put {$#2$} [b] <0pt, 5pt> at -2 3
\put {$#3$} [tl] <3pt, 0pt> at -2.5 3.5
\put {$#4$} [b] <0pt, 5pt> at -1 2
\put {$#5$} [tl] <3pt, 0pt> at -1.5 2.5
\put {$#6$} [b] <0pt, 5pt> at 0 1
\put {$#7$} [tl] <3pt, 0pt> at -0.5 1.5
\put {$#8$} [Br] <-3pt, 0pt> at 0 -0.5
\endpicture}

\def\ctreec#1#2#3#4#5#6#7{
\beginpicture
\small
\setcoordinatesystem units <0.5cm,0.6cm>
\linethickness=1pt
\setplotsymbol ({\tenrm .})
\setlinear
\plot 0 -0.5  0 1  -3 4 /
\plot 0 1  3 4 /
\plot -2 3  -1 4 /
\plot 2 3  1 4 /
\put {$#1$} [b] <0pt, 3pt> at -3 4
\put {$#1$} [b] <0pt, 3pt> at -1 4
\put {$#1$} [b] <0pt, 3pt> at 1 4
\put {$#1$} [b] <0pt, 3pt> at 3 4
\put {$#2$} [b] <0pt, 5pt> at -2 3
\put {$#3$} [b] <0pt, 5pt> at 2 3
\put {$#4$} [tr] <-3pt, 0pt> at -1 2
\put {$#5$} [tl] <3pt, 0pt> at 1 2
\put {$#6$} [b] <0pt, 5pt> at 0 1
\put {$#7$} [Br] <-3pt, 0pt> at 0 -0.5
\endpicture}

\def\ltreed#1#2#3#4#5{
\beginpicture
\small
\setcoordinatesystem units <0.5cm,0.6cm>
\linethickness=1pt
\setplotsymbol ({\tenrm .})
\setlinear
\plot 0 -0.5  0 1  -3 4 /
\plot 0 1  3 4 /
\plot -1 2  1 4 /
\plot 0 3  -1 4 /
\put {$#1$} [b] <0pt, 3pt> at -3 4
\put {$#1$} [b] <0pt, 3pt> at -1 4
\put {$#1$} [b] <0pt, 3pt> at 1 4
\put {$#1$} [b] <0pt, 3pt> at 3 4
\put {$#2$} [b] <0pt, 5pt> at 0 3
\put {$#3$} [tl] <0pt, 0pt> at -0.5 2.5
\put {$#4$} [tr] <-3pt, 0pt> at -0.5 1.5
\put {$#5$} [Br] <-3pt, 0pt> at 0 -0.5
\endpicture}

\def\rtreed#1#2#3#4#5{
\beginpicture
\small
\setcoordinatesystem units <-0.5cm,0.6cm>
\linethickness=1pt
\setplotsymbol ({\tenrm .})
\setlinear
\plot 0 -0.5  0 1  -3 4 /
\plot 0 1  3 4 /
\plot -1 2  1 4 /
\plot 0 3  -1 4 /
\put {$#1$} [b] <0pt, 3pt> at -3 4
\put {$#1$} [b] <0pt, 3pt> at -1 4
\put {$#1$} [b] <0pt, 3pt> at 1 4
\put {$#1$} [b] <0pt, 3pt> at 3 4
\put {$#2$} [b] <0pt, 4pt> at 0 3
\put {$#3$} [tr] <0pt, 0pt> at -0.5 2.5
\put {$#4$} [tl] <3pt, 0pt> at -0.5 1.5
\put {$#5$} [Br] <-3pt, 0pt> at 0 -0.5
\endpicture}



\def\implies{\Rightarrow}
\def\ideq{\equiv}
\def\hom{\hbox{hom}}

\def\inv{^{\scriptscriptstyle-\!1}}

\def\Z{\Bbb Z}
\def\G{{\cal G}}

\def\tens{\otimes}
\def\id{\hbox{\rm ID}}
\def\ci{^{\circ\hbox{--}1}}
\def\re{r_\epsilon}
\def\ce{c_\epsilon}
\def\lam{\lambda(\epsilon)}
\def\bas#1#2#3{(#1\!\!\to\!\!#2#3)}
\def\ato{\smash{\mathop{\longmapsto}\limits^{\alpha}}\ }
\def\mee{\mu[\epsilon,\epsilon]}
\def\labmapsto#1{\smash{\mathop{\hbox to 30pt{\rightarrowfill}}\limits^{#1}}}
\def\field{F_{p^\alpha}}
\def\fadd{F_{p^\alpha}^+}
\def\fmul{F_{p^\alpha}^*}

\def\ltimes{\triangleright\!\!\!<}
\def\iso{\cong}
\def\tens{\otimes}
\def\unit{\hbox{$\epsilon$}}
\def\dsum{\oplus}
\def\cite#1{\rm{[#1]}}
\def \pit{+\hbox {irrelevant terms}}


\title{Near-group categories}
\author{Jacob Siehler}
\address{Department of Mathematics, Virginia Tech\\Blacksburg, 
VA 24061-0123, USA} 
\email{jsiehler@math.vt.edu}

\abstract 
We consider the possibility of semisimple tensor categories whose fusion
rule includes exactly one noninvertible simple object.  Conditions are
given for the existence or nonexistence of coherent associative structures
for such fusion rules, and an explicit construction of matrix solutions to
the pentagon equations in the cases where we establish existence.  Many of
these also support (braided) commutative and tortile structures and we
indicate when this is possible.  Small examples are presented in detail.
\endabstract

\primaryclass{18D10}

\maketitle

\section{Introduction and results}

The term {\it near-group categories} is introduced to describe a specific
class of finite semisimple monoidal categories with duality.  For the
purposes of the present paper:

{\it Finite semisimple} means: our categories have a ``ground ring'' $R$,
and hom sets in the category are free $R$--modules;  that is, there is
a finite set of simple objects ${s_i}$ with $$\hom(s_i,s_j)\iso\cases
{R,&$i=j$\cr 0,& $i\ne j$} $$ and every object in the category is
(isomorphic to) a direct sum of simple objects.

{\it Monoidal} means: there is a (bifunctorial) tensor product $\tens\co\cal
C \times \cal C \to \cal C$.  We assume the existence of a unit object
(\unit) for our tensor product (so that tensoring on either side by
\unit\ is naturally isomorphic to the identity).  For objects $a,b,c$
there are natural associativity morphisms $\alpha_{a,b,c}\co (a\tens b)\tens
c\to a\tens(b\tens c)$, and these satisfy the pentagon axiom.

{\it Duality} means:  every object $x$ has a dual object $x^*$.  Moreover,
there exist ``pairing'' and ``copairing'' morphisms $\lambda_x\co x^*\tens
x\to \unit$ and $\Lambda_x\co \unit\to x\tens x^*$ which satisfy the
identities:
\eject
$$\left(\id_x\tens\lambda_x\right)\left(\alpha_{x,x^*,x}\right)
  \left(\Lambda_x\tens\id_x\right)=\id_x$$
and
$$\left(\lambda_x\tens\id_{x^*}\right)\left(\alpha_{x^*,x,x^*}\right)
  \inv\left(\id_{x^*}\tens\Lambda_x\right)=\id_{x^*}$$
A consequence of this hypothesis is the existence of adjunction isomorphisms
$\hom(x\tens y,z)\iso \hom(x,z\tens y^*)$ and $\hom(x,y\tens z)\iso
\hom(y^*\tens x,z)$.  As a special case, simples $s_i$ and $s_j$ have
$$\hom(\unit,s_i\tens s_j)\iso\cases{R,&$s_j=s_i^*$\cr 0,&otherwise} $$
Categories in which every simple has a multiplicative inverse (so
$g\tens g\inv \iso \unit$) have been studied carefully, and called
variously {\it $\Theta$--categories} \cite {FK} or {\it group--categories}
\cite {Q2}.  Note that the multiplicative inverse is the same object as
the dual. Complete information is available about the classification
of categorical structures on group--categories and the interpretation
of their field theories (in \cite {Q2}). The classification is possible
in one sense because all the structural equations are 1--dimensional and
there is no non-commutative matrix arithmetic involved in their solutions.

In the present paper we would like to advance to studying a slightly
more complicated sort of multiplication, in which there is (up to
isomorphism) a single noninvertible simple object.  For a category with
such a multiplication we introduce the term {\it near-group category}.
The noninvertible simple object we will call $m$.  The invertible simple
objects comprise a group $G$ under tensor product; invertibility of the
group elements implies
$$g \tens m \iso m \tens g \iso m$$
for any invertible $g$.  Moreover, the duality assumption implies 
$$m\tens m \iso G \dsum k\cdot m$$
that is, $m\tens m$ splits as one summand of each invertible type,
and some number $k$ of $m$--summands (possibly $k=0$).  We see that the
multiplication table (or {\it fusion rule}) for a near-group category
is described by giving the group $G$ and the integer $k$; we will refer
to categories with the {\it near-group fusion rule $(G,k)$}.

The first question of interest is which $(G,k)$ actually occur
as fusion rules of near-group categories; the primary issue is the
existence of coherent associativity, but we also want to investigate
which $G$ and $k$ admit the additional layers of structure we have 
indicated in the preceding section.  The case where $k=0$ is dealt
with in \cite{TY} and \cite{S}, so the present paper is concerned
only with near-group fusion rules having $k\ge 1$.

The main results on monoidal structures are:

\medskip
{\bf 1.1\qua Theorem}\qua{\rm(Order control)}\qua\sl
If the near-group fusion rule $(G,k)$ admits a monoidal structure then
$|G|\le k+1$.\rm

\medskip
{\bf 1.2\qua Theorem}\qua{\rm(Existence)}\qua\sl
In the maximal case $|G|=k+1$ the near-group fusion rule $(G,k)$ admits a
monoidal structure iff $G$ is the multiplicative group of a finite field
(ie, cyclic of order $p^\alpha-1$).\rm

\medskip
To prove necessity, we examine in detail (sections 4 and 5) the pentagon
equations in the category; we find that $G$ must support some special
structure (4.2) and we show (6.1) that only the specified groups support
this structure.


We give two different proofs of sufficiency. In 4.1 a finite
group is constructed whose category of representations has the given 
fusion rule. Representation categories are symmetric monoidal, so
this also shows if there is any monoidal category then there is a 
symmetric one. The other proof of sufficiency is in 5.18. This uses
the analysis of the pentagon equations carried out in sections 4 and 
5 to directly exhibit solutions of the matrix equations. This
approach gives information about the whole set of categories, 
including some that do not admit braidings. There is a ``simplest'' 
standard solution; analysis of hexagon equations in sections 8 and 
9 verifies that this solution additionally admits a symmetric commuting
structure.  It seems likely that this is the representation category 
exhibited earlier, but this has not been verified.

The analysis of the pentagon equations enables small cases to be 
worked out completely.
\proclaim{1.3 Theorem}
The first three maximal-order cases have the following structure:
\bigskip
\def\medstrut{\vrule width 0pt height 10pt depth 5pt}
\def\lgstrut{\vrule width 0pt height 21pt depth 5pt}
\def\hgstrut{\vrule width 0pt height 35pt depth 5pt}
\vbox{\offinterlineskip\small\sl
\halign{\vrule#&~#~&\vrule#&~#&\vrule#&~#&\vrule#&~#&\vrule#\cr
\noalign{\hrule}
& \medstrut Fusion\hfil && $(\Bbb Z/2,1)$ && $(\Bbb Z/3,2)$ && $(\Bbb Z/4,3)$ &\cr
\noalign{\hrule}
& \medstrut Field\hfil  && $F_3$          && $F_{2^2}$      && $F_5$          &\cr
\noalign{\hrule}
& \lgstrut\vbox{\halign{#\cr
                Monoidal\hfil\cr
                structures\hfil\cr}}\hfil &&
  \vbox{\halign{#\cr
  3, indexed by $\xi$\cr
  where $\xi^3 = 1$\hfil\cr}} &&
  \vbox{\halign{#\cr
  2, indexed by $\xi = \pm 1$\cr
  \hfil\cr}} &&
  \vbox{\halign{#\cr
  unique monoidal\hfil\cr
  structure\cr}} &\cr
\noalign{\hrule}
& \hgstrut\vbox{\halign{\strut#\cr
                Braidings\hfil\cr
                \hfil\cr
                \hfil\cr}}\hfil && 
  \vbox{\halign{#\cr
                $\xi = 1$: 3, indexed\cr
                \hfil by $\psi^3 = 1$\cr
                $\xi \not= 1$: not braided\cr}} &&
  \vbox{\halign{#\cr
                $\xi = 1$: 4, indexed\cr
                \hfil by $\psi(1),\psi(2) = \pm 1$\cr
                $\xi \not= 1$: not braided\cr}} &&
  \vbox{\halign{#\cr
                unique\hfil\cr
                braiding\hfil\cr
                \hfil\cr}} &\cr
\noalign{\hrule}
& \lgstrut\vbox{\halign{#\cr
                Balance\hfil\cr
                \hfil\cr}}\hfil &&
  \vbox{\halign{#\cr
                $\psi = 1$: balanced\hfil\cr
                $\psi \not= 1$: not balanced\hfil\cr}} &&
  \vbox{\halign{#\cr
                all balanced\hfil\cr
                \hfil\cr}} &&
  \vbox{\halign{#\cr
                balanced\hfil\cr
                \hfil\cr}} &\cr
\noalign{\hrule}
& \lgstrut\vbox{\halign{#\cr
        Symmetry\hfil\cr
        \hfil\cr}}\hfil &&
  \vbox{\halign{#\cr
                $\psi = 1$: symmetric\hfil\cr
                \hfil\cr}} &&
  \vbox{\halign{#\cr
                $\psi(1) = \psi(2)$: symmetric\hfil\cr
                $\psi(1) = -\psi(2)$: not symmetric\cr}} &&
  \vbox{\halign{#\cr
                symmetric\hfil\cr
                \hfil\cr}} &\cr
\noalign{\hrule}
}}\hfil
\medskip
\endproclaim
The primes 2 and 3 play special roles in the analysis, so the 
complexity of the first two cases probably is representative of the 
fusion rules corresponding to finite fields of characteristic 2 or 3. 
Larger primes may have simpler behavior, though the complete uniqueness 
seen in the $F_5$ case is probably too simple.

At the other extreme from the maximal-order case we have:

\proclaim {1.4 Theorem}
Assume that the characteristic of the ground ring $R$ is not equal to 2.
Then, in the minimal case $|G|=1$, the near-group fusion rule $(G,k)$
does not admit a monoidal structure if
$k\equiv 2 \hbox{\ or\ }3\hbox{\rm\ mod\ }4$
\endproclaim

\proclaim {Remark} {\rm A stronger result is obtained in \cite{O} for
the case where the ground ring has characteristic zero, but the present
proof works also for positive odd characteristic.}
\endproclaim

We prove this in section 7.  When $G$ is small there is less symmetry
to exploit and the pentagon equations become much more difficult;
nonexistence theorems are easier than constructions.  For $k=1$ and
$|G|=1$ it is not difficult to solve (the unique characteristic 5
solution is given in \cite Q ; this lifts to two distinct solutions
in characteristic 0). But the next possible trivial-group case, $k=4$,
so far seems to be intractable.

Section 8 reduces the hexagon equations and gives a standard construction
of commutativity data.  Section 9 then gives the braidings possible for
the example categories from section 3, and this includes nonsymmetric
braidings on the standard monoidal structures; it turns out that the
more exotic monoidal structures exhibited in section 3 do not support
braidings.

Section 10 gives a brief note on the possibility of adding twist
morphisms to the braided near-group categories studied in sections 8 and 9
(producing ``tortile categories'').

\section{Notation for associativities}

In section 3 we will give examples and complete data for categories with
small values of $k$.  Before doing that we need to introduce notation used
to present the many different associativitities in the category; this
scheme is an extension of the notation used in \cite{TY} , and in later
sections this notation is used extensively in proofs and calculations.
Throughout this section, let $a,b,c$ stand for arbitrary elements of $G$.

\sh{2.1\qua Associativities involving a product of three group elements}
These we denote by a function $\alpha$ depending on three group variables. 
$$\ltreea abc{}{ab}{}{abc}\qua
  \raise 20pt\hbox{$\labmapsto {\alpha(a,b,c)}$}\qua
  \rtreea abc{}{bc}{}{abc}
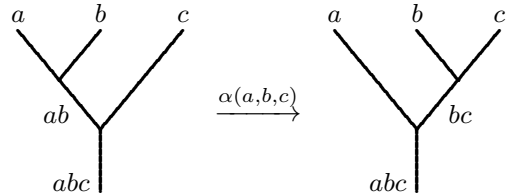
$$

\sh {2.2\qua Associativities involving two group elements and one $m$}
These we think of as a function of the two group elements, denoted
by $\alpha$'s with a subscript to indiciate the position of the $m$:
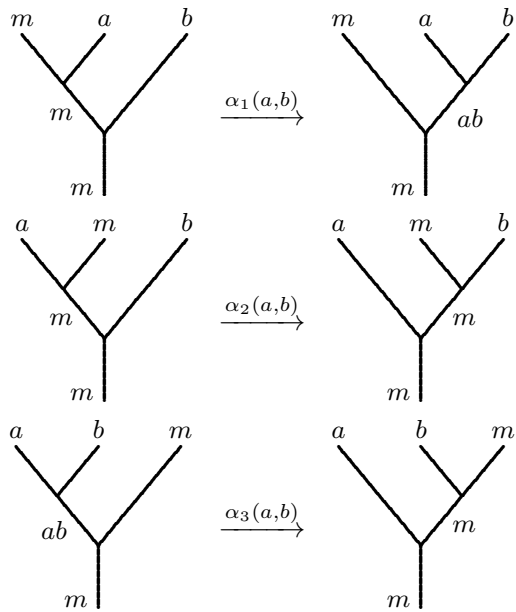
$$\ltreea mab{}m{}m\qua
  \raise 20pt \hbox{$\labmapsto{\alpha_1(a,b)}$}\qua
  \rtreea mab{}{ab}{}m$$
$$\ltreea amb{}m{}m\qua
  \raise 20pt \hbox{$\labmapsto{\alpha_2(a,b)}$}\qua
  \rtreea amb{}m{}m$$
$$\ltreea abm{}{ab}{}m\qua
  \raise 20pt \hbox{$\labmapsto{\alpha_3(a,b)}$}\qua
  \rtreea abm{}m{}m$$

\sh {2.3\qua Associativities involving two $m$'s and one group element} 
Such a product has both group summands and $m$ summands in it.  We will
use $\beta$'s to denote the associativity on the group summands:
$$\ltreea amm{}m{}b\qua
  \raise 20pt\hbox{$\labmapsto {\beta_1(a,b)}$}\qua
  \rtreea amm{}{a^{-1}b}{}b
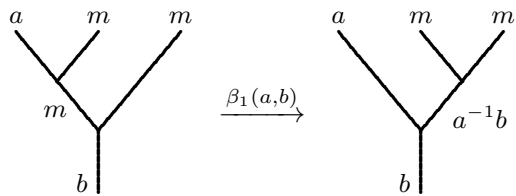
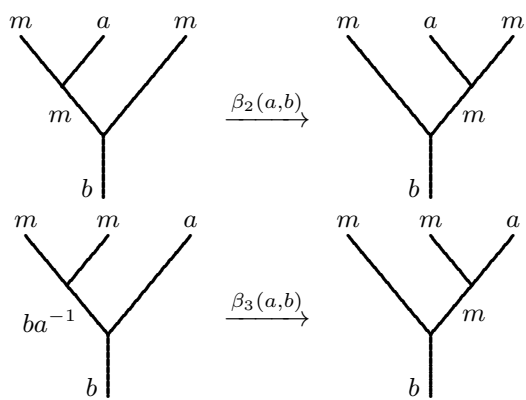
$$
$$\ltreea mam{}m{}b\qua
  \raise 20pt\hbox{$\labmapsto {\beta_2(a,b)}$}\qua
  \rtreea mam{}m{}b$$
$$\ltreea mma{}{ba^{-1}}{}b\qua
  \raise 20pt\hbox{$\labmapsto {\beta_3(a,b)}$}\qua
  \rtreea mma{}m{}b$$
For the same kind of product, $\gamma$'s will denote the the associativity
on the $m$ summands.  These will take values in $k\times k$ matrices
(rows indexed by $i=1\ldots k$, columns indexed by $j=1\ldots k$):
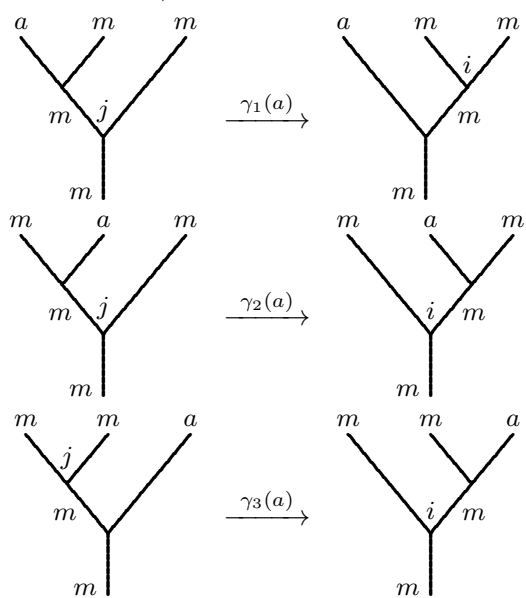
$$\ltreea amm{}mjm\qua
  \raise 20pt\hbox{$\labmapsto {\gamma_1(a)}$}\qua
  \rtreea ammim{}m$$
$$\ltreea mam{}mjm\qua
  \raise 20pt\hbox{$\labmapsto {\gamma_2(a)}$}\qua
  \rtreea mam{}mim$$
$$\ltreea mmajm{}m\qua
  \raise 20pt\hbox{$\labmapsto {\gamma_3(a)}$}\qua
  \rtreea mma{}mim$$

\sh {2.4\qua Associativities for a threefold product of $m$'s} 
Such a product has both group summands and $m$ summands.  We will use
$\lambda$ to denote associativity on the group summands.  The $\lambda$'s
will be $k\times k$ matrices:
$$\ltreea mmmjm{}a\qua
  \raise 20pt\hbox{$\labmapsto {\lambda(a)}$}\qua
  \rtreea mmmim{}a$$
Finally, $\mu$ stands for associativity on the $m$ summands; $\mu$ will be
a $(k^2+|G|)\times(k^2+|G|)$ matrix:  on the left, there are $|G|$
basis elements as $x$ runs over $G$, and $k^2$ indexed by $r,s=1\ldots k$
when $x=m$.  Similarly on the right:
$$\ltreea mmmsxrm\qua
  \raise 20pt\hbox{$\labmapsto {\mu}$}\qua
  \rtreea mmmjyim$$

\section{Examples and results for small $k$}

In this section we would like to present, for certain groups $G$ and
small values of $k$, the complete classification of monoidal categories
with fusion rule $(G,k)$ by giving explicit associativity matrices.
The purpose is to make numerical data conveniently available in the
small cases and to illustrate the construction given in section 5.18.
The matrix construction is best explained using a permutation $\pi$
which acts on the nonidentity elements of the group; this satisfies some
characterizing identities which are set out in 4.2 but are not important
for the present purpose.

\sh{3.1\qua Example 1\qua $k=1, G={\Bbb Z}/2{\Bbb Z}=\{\epsilon,g\}$ }

Here there is only one nonidentity element in $G$ so of course $\pi$
is the trivial permutation.

In fact, there are three distinct monoidal structures possible for this
fusion rule, classified by a choice of $\xi$ with $\xi^3=1$ (of course
in characteristic 3 these collapse to a single solution).

Let $\chi$ be the nontrivial character of $G$, that is,
$\chi(\epsilon)=1$, $\chi(g)=-1$.

Let $a,b,c$ represent arbitrary elements of $G$. With a good choice of
basis, the associativities are as follows:
$$\eqalign{
\alpha(a,b,c)&\ideq\pmatrix{1}\cr
\alpha_1(a,b)=\alpha_2(a,b)=\alpha_3(a,b)&\ideq\pmatrix{1}\cr
\beta_1(a,b)=\beta_2(a,b)=\beta_3(a,b)&\ideq\pmatrix{1}\cr
\gamma_1(a)=\gamma_2(a)=\gamma_3(a)&=\pmatrix{\chi(a)}\cr
\lambda(a) &= \pmatrix{\xi\chi(a)}\cr}$$
\noindent and finally, $\mu$:
$$\bordermatrix{
\hfill x=&\epsilon&g&m \cr
y = \epsilon &{1\over 2}&{1\over 2}&{1\over 2\xi} \cr
\hfill g&{1\over 2}&{1\over 2}&-{1\over 2\xi} \cr
\hfill m&1&-1&0\cr}$$
{\bf Note}\qua The construction in 5.18 corresponds to $\xi=1$, which is also
the structure of the category of representations of the group $S_3$,
in characteristic prime to 6.

\sh{3.2\qua Example 2\qua $k=2, G={\Bbb Z}/3{\Bbb Z} = \{\epsilon,g,g^2\}$}

In this case take $\pi$ to be the identity permutation on $\{g,g^2\}$.  

Let $\chi_1$ and $\chi_2$ be the two nontrivial characters of $G$.

There are two distinct monoidal structures for this fusion rule,
corresponding to a choice of $\xi=\pm1$, and with a good choice of bases
the associativity data is as follows:
$$\eqalign{
\alpha(a,b,c)&\ideq\pmatrix{1}\cr
\alpha_1(a,b)=\alpha_2(a,b)=\alpha_3(a,b)&\ideq\pmatrix{1}\cr
\beta_1(a,b)=\beta_2(a,b)=\beta_3(a,b)&\ideq\pmatrix{1}\cr
\gamma_1(a)=\gamma_2(a\inv)=\gamma_3(a)&=\pmatrix{\chi_1(a)&  0 \cr 0 & \chi_2(a)\cr} \cr
\lambda(a)&=\pmatrix{\xi\chi_1(a) & 0 \cr 0 & \xi\chi_2(a) }\cr}$$
and $\mu$:

$$\bordermatrix{
&\smallstrut\smash{\raise10pt\rlap{$x=$\hss}}\epsilon&g&g^2&\smash{\raise12pt
\rlap{$(r,s)=$\hss}(1,1)}&(1,2)&(2,1)&(2,2)\cr
\hfill y=\epsilon&\bigstrut{\displaystyle{\xi\over3}}&{\displaystyle{
  \xi\over3}}&{\displaystyle{ \xi\over3}}&0&\displaystyle{
  \xi\over3}&{\displaystyle{1\over3}}&0\cr
\hfill g&\bigstrut{\displaystyle{ \xi\over3}}&{\displaystyle{
  \xi\over3}}&{\displaystyle{ \xi\over3}}&0&\displaystyle{
  \xi\chi_1(g^2)\over3}&\displaystyle{\chi_2(g^2)\over3}&0\cr
\hfill g^2&\bigstrut{\displaystyle{ \xi\over3}}&{\displaystyle{
  \xi\over3}}&{\displaystyle{ \xi\over3}}&0&\displaystyle{
  \xi\chi_1(g)\over3}&\displaystyle{\chi_2(g)\over3}&0\cr
(i,j)=(1,1)&\smallstrut0&0&0&0&0&0&\xi\cr
\hfill(1,2)&\smallstrut1&\chi_1(g^2)&\chi_1(g)&0&0&0&0\cr
\hfill(2,1)&\smallstrut\xi&\xi\chi_2(g^2)&\xi\chi_2(g)&0&0&0&0\cr
\hfill(2,2)&\smallstrut0&0&0&1&0&0&0\cr}$$
{\bf Note}\qua The construction in 4.18 corresponds to $\xi=+1$, and this is
also the structure in the category of representations of the group $A_4$
in characteristic prime to 12.

\sh{3.3\qua Example 3\qua $k=3, G={\Bbb Z}/4{\Bbb Z} = \{\epsilon,g,g^2,g^3\}$ }

The data for this example is lengthy, but this is the smallest available
example where $\pi$ is a nontrivial permutation, and is really the best
illustration of the construction in 5.18.

Take $\pi$ to be the 3--cycle $\left(g\ g^1\ g^3\right )$.

Let $\chi_1,\chi_2,$ and $\chi_3$ be the nontrivial characters of $G$,
with $\chi_2(g^2)=1$.

There is a unique monoidal structure for this fusion rule, and with a good
choice of bases the associativity data is as follows:
$$\eqalign{
\alpha(a,b,c)&\ideq\pmatrix{1}\cr
\alpha_1(a,b)=\alpha_2(a,b)=\alpha_3(a,b)&\ideq\pmatrix {1}\cr
\beta_1(a,b)=\beta_2(a,b)=\beta_3(a,b)&\ideq\pmatrix {1}\cr
\gamma_1(a)&=\pmatrix {\chi_1(a)&0&0 \cr 0&\chi_2(a)&0 \cr 0&0&\chi_3(a)\cr}\cr
\lambda(\epsilon)&=\pmatrix {0&0&1 \cr 1&0&0 \cr 0&1&0 \cr}\cr
\gamma_2(a)&=\lambda(\epsilon)\inv \gamma_1(a\inv) \lambda(\epsilon)\cr
\gamma_3(a)&=\lambda(\epsilon)\inv \gamma_2(a\inv) \lambda(\epsilon)\cr
\lambda(a)&=\lambda(\epsilon)\gamma_1(a)\cr}$$
and the large associator $\mu$ is given in four submatrices: the upper
left ``$M$ submatrix'' corresponding to $x,y\in G$:
\def\tinystrut{\vrule width 0pt height 5pt depth 5pt}
\def\qu{\tinystrut {1\over4}}
$$\bordermatrix {
\hfill x=&\epsilon&g&g^2&g^3 \cr
y=\epsilon&\qu&\qu&\qu&\qu\cr
\hfill g&\qu&\qu&\qu&\qu\cr
\hfill g^2&\qu&\qu&\qu&\qu\cr
\hfill g^3&\qu&\qu&\qu&\qu\cr}$$
the upper right ``$R$ submatrix'', corresponding to $x=m, y\in G$:
$$
{\small\bordermatrix {
\hfill(r,s)= & (1,1)\!\!&(1,2)\!\! & (1,3) & (2,1)\!\! & (2,2)\!\! & (2,3) &
  (3,1)\!\! & (3,2)\!\! & (3,3) \cr
\hfill y=\epsilon & \smallstrut 0 & 0 & \qu & 0 & \qu & 0 & \qu & 0 & 0 \cr
\hfill g &\smallstrut 0 & 0 & \llap{$1\over4$}\chi_3(g^3) & 0 &
  \llap{$1\over4$}\chi_1(g^3) & 0 & \llap{$1\over4$}\chi_2(g^3) & 0 & 0 \cr
\hfill g^2 &\smallstrut 0 & 0 & \llap{$1\over4$}\chi_3(g^2) & 0 &
  \llap{$1\over4$}\chi_1(g^2) & 0 & \llap{$1\over4$}\chi_2(g^2) & 0 & 0 \cr
\hfill g^3 &\smallstrut 0 & 0 & \llap{$1\over4$}\chi_3(g) & 0 &
  \llap{$1\over4$}\chi_1(g) & 0 & \llap{$1\over4$}\chi_2(g) & 0 & 0 \cr
}}$$
the lower left ``$C$ submatrix'', corresponding to $x\in G, y=m$:
$$\bordermatrix {
\hfill x=&\epsilon&g&g^2&g^3 \cr
(i,j)=(1,1)&0&0&0&0\cr
\hfill(1,2)&1&\chi_1(g^3)&\chi_1(g^2)&\chi_1(g)\cr
\hfill(1,3)&0&0&0&0\cr
\cr
\hfill(2,1)&1&\chi_2(g^3)&\chi_2(g^2)&\chi_2(g)\cr
\hfill(2,2)&0&0&0&0\cr
\hfill(2,3)&0&0&0&0\cr
\cr
\hfill(3,1)&0&0&0&0\cr
\hfill(3,2)&0&0&0&0\cr
\hfill(3,3)&1&\chi_3(g^3)&\chi_3(g^2)&\chi_3(g)\cr}$$
and the lower right ``$N$ submatrix'' corresponding to $x=y=m$, with
columns indexed by $(r,s)$ and rows indexed by $(i,j)$ in the same
(lexicographic) order they appear in the previous two submatrices:
$$\pmatrix{ 
0&0&0\quad  &0&0&0\quad  &0&1&0 \cr
0&0&0\quad  &0&0&0\quad  &0&0&0 \cr
0&0&0\quad  &0&0&1\quad  &0&0&0 \cr
\cr
0&0&0\quad  &0&0&0\quad  &0&0&0 \cr
0&0&0\quad  &0&0&0\quad  &0&0&1 \cr
1&0&0\quad  &0&0&0\quad  &0&0&0 \cr
\cr
0&0&0\quad  &1&0&0\quad  &0&0&0\strut \cr
0&1&0\quad  &0&0&0\quad  &0&0&0 \cr
0&0&0\quad  &0&0&0\quad  &0&0&0 \cr}$$
This example illustrates the most important structural features of the
solutions in the maximal group case, ie, in the setting of theorem 1.2:
$\alpha$ and $\beta$ associativities are trivial; the $\gamma$'s are
$k$--dimensional representations of $G$ built out of the nontrivial
characters of $G$; the shape of the $\lambda(\epsilon)$ matrix comes
from the permutation $\pi$.

Moreover, the large associator $\mu$ is predictably organized.
The $M$ submatrix is constant.  The $C$ submatrix has one nonzero row
in each block (ie, one nonzero row for each $i=1,2,\ldots,k$), and the
columns are simply related by characters of $G$.  The $R$ submatrix is
organized much like $C$.

Finally, if we identify $G$ with its dual group of characters,
$\pi$ can be viewed as acting on $\chi_1,\chi_2,\chi_3$ as the
3--cycle $\left(\chi_1\ \chi_2\ \chi_3\right)$.  This is the key
to understanding the $N$ submatrix: a nonzero entry occurs in row
$(i,j)$, column $(r,s)$ if and only if $\chi_r\chi_s=\chi_i$ and
$\pi\inv(\chi_i)\pi\inv(\chi_j)=\pi\inv(\chi_r)$.

Elaborating these observations and turning them into theorems will be
the work of the next two sections.

\section{Pentagon equations and proof of Theorem 1.2}

\sh{4.1\qua Proof of theorem 1.2: sufficiency}
For each finite field $\field$, we construct a nonabelian group whose
category of representations is a  near-group category in which the group
of invertibles is isomorphic to the multiplicative group of the field.

Write $\fadd$ for the additive group of the field with $p^\alpha$
elements, and $\fmul$ for the multiplicative group.  Let $\cal G$ be
the semidirect product $\fadd\ltimes\fmul$ with the obvious action.

It's easy to verify:

(1)\qua  $[\G,\G]=\fadd$

(2)\qua  $G$ has one conjugacy class containing the identity, one
conjugacy class consisting of the nonidentity elements in $\fadd$,
and $(p^\alpha-2)$ classes each of order $p^\alpha$.

It follows that in the semisimple setting, $\G$ has $(p^\alpha-1)$
linear representations comprising a cyclic group $\Z_{p^\alpha-1}$,
together with a single, noninvertible, $(p^\alpha-1)$--dimensional
representation $m$.  By dimension, $m\tens m$ contains $(p^\alpha-2)$
copies of $m$, so the category of representations of $\G$ has near-group
fusion rule $\left(\Z_{p^\alpha-1},p^\alpha-2\right)$, and the ``if''
of Theorem 1.2 follows.\endprf

\proclaim {Remark} {\rm The group $\G$ is none other than the affine
group of the field $\field$.  The presentation as a semidirect product
is just convenient for counting its representations.}
\endproclaim

We will later (5.18) show explicitly how to build associativity matrices
that solve the pentagon equations for these fusion rules.  That point of
view emerges from the analysis of the pentagon equations we carry out
in this and the following section.  Our goal is to prove the following
intermediate theorem:

\proclaim {4.2\qua Intermediate Theorem}  
If the near-group fusion rule $(G,k)$, $|G|=k+1$, supports a monoidal
structure, then there is a permutation $\pi$ on the nonidentity elements
of $G$ with the following three properties:

{\rm (i)}\qua $\pi^3=\hbox{id}$

{\rm (ii)}\qua $\pi(x)\inv = \pi\inv(x\inv)$

{\rm (iii)}\qua $\pi(st)=\pi(t)\pi\bigl[\pi(s)\inv\pi(t\inv)]$ (for all
$s\ne t\inv$)
\endproclaim

Once we have established the intermediate theorem, we analyze (6.1)
which groups $G$ support the structure of such a permutation $\pi$,
completing the ``only if'' portion of Theorem 1.2.  It is interesting
to note that the factorization property (iii) on the permutation $\pi$
is similar to the property studied by \cite{KR} (particularly, their
Proposition 2) where they produce solutions to the pentagon equations from
``symmetrically factorizable groups.''

\sh{ Pentagon equations}
Assume we have a near-group category with fusion rule $(G,k)$.  

Equations (1)--(11) of \cite{TY} hold verbatim in this setting;
consequently after fixing bases for all hom sets except $\hom(m,mm)$
we may assume:

$\alpha\ideq\alpha_1\ideq\alpha_3\ideq1$

$\beta_1\ideq1$

$\beta_2\ideq\alpha_2$ are symmetric and bimultiplicative

$\beta_3(a,b)=\beta_3(a,\epsilon)$

Since $\beta_3$ is independent of the second variable, we will abbreviate
$\beta_3(a,b)$ to simply $\beta_3(a)$.

We proceed by examining remaining pentagon equations and looking for a
good choice of basis for $\hom(m,mm)$.  In what follows we will refer to
``the pentagon $abcd/x$'', meaning the content of the pentagon equation
for the $x$ summands in the product $a\tens b\tens c\tens d$.

\proclaim {4.3\qua Proposition}
$\gamma_1$, $\gamma_2$, and $\gamma_3$ are representations of $G$.
\endproclaim

\prf
For $a,b\in G$, this is precisely the content of the pentagons $mmab/m$,
$abmm/m$, and $mabm/m$, respectively:
$$\eqalignno{
\gamma_3(b)\gamma_3(a)&=\gamma_3(ab)&(1)\cr
\gamma_1(b)\gamma_1(a)&=\gamma_1(ab)&(2)\cr
\gamma_2(b)\gamma_2(a)&=\gamma_2(ab)&(3)\cr}$$
Pentagons $ammb/m$, $mamb/m$, and $ambm/m$ are commutator relationships among the $\gamma$ representations:
$$\eqalignno{
\gamma_3(b)\alpha_2(a,b)\gamma_1(a)&=\gamma_1(a)\gamma_3(b)&(4)\cr
\alpha_2(a,b)\gamma_3(b)\gamma_2(a)&=\gamma_2(a)\gamma_3(b)&(5)\cr
\gamma_2(b)\gamma_1(a)\alpha_2(a,b)&=\gamma_1(a)\gamma_2(b)&(6)\cr}$$
\endprf

\proclaim {4.4\qua Proposition}
$\alpha_2(a,b)$ is an $r$--th root of unity where $r|\gcd(k,|G|)$.
\endproclaim

\prf
We choose basis of $\hom(m,mm)$ to split the representation
$\gamma_1$ into $k$ linear summands, diagonalizing
all $\gamma_1$ matrices.  Now, equation $(6)$ says
$$\gamma_2(b)\inv\gamma_1(a)\gamma_2(b)=\alpha_2(a,b)\gamma_1(a)$$
The right hand side is of course diagonal; the left hand side is therefore
diagonal and has the same eigenvalues as $\gamma_1(a)$.  Multiplication by
$\alpha_2(a,b)$ simply permutes those eigenvalues.

Since $\alpha_2$ is multiplicative in each factor we know $\alpha_2(a,b)$
is an $r$--th root of unity for some $r$ which divides $|G|$.  But by the
above observation the $k$ eigenvalues of $\gamma_1(a)$ break up into
$n$ orbits each of size $r$;  we obtain $k=nr$ and so $r$ divides $k$
as well.
\endprf

\proclaim {4.5\qua Proposition}
The representations $\gamma_1$, $\gamma_2\inv$, and $\gamma_3$ are all
conjugate, and $\lambda(\epsilon)$ is an intertwiner carrying $\gamma_3$
to $\gamma_1$ to $\gamma_2\inv$ back to $\gamma_3$.
\endproclaim

\prf
The pentagons for $mmmg/hg$, $gmmm/hg$, $mmgm/h$, and $mgmm/h$ give the
following equations respectively:
$$\eqalignno{
\lambda(hg)&=\gamma_3(g)\lambda(h)&(7)\cr
\lambda(hg)&=\lambda(h)\gamma_1(g)&(8)\cr
\gamma_2(g)\lambda(h)\gamma_3(g)&=\lambda(h)\beta_2(g,h)&(9)\cr
\gamma_1(g)\lambda(h)\gamma_2(g)&=\lambda(h)\beta_2(g,h)&(10)\cr
}$$
By (7) with $h=\epsilon$, for any $g$ we have
$\lambda(g)=\gamma_3(g)\lambda(\epsilon)$.

Now substitute into (8) with $h=\epsilon$ again to get
$\gamma_1(g)=\lambda(\epsilon)\inv\gamma_3(g)\lambda(\epsilon)$, so
$\lambda(\epsilon)$ carries $\gamma_3$ to $\gamma_1$. Take $h=\epsilon$
in the remaining equations to complete the remaining claims.
\endprf

\sh{Notation for the associator $\mu$}

Recall from section 2.4 that $\mu$ is the associator for the threefold
product of $m$'s, on its $m$ summands.  We think of it concretely as a
matrix with respect to the standard bases indicated by the trees in 2.4.
On the left, there are $|G|$ basis elements indexed by the elements of
$G$, and $k^2$ elements indexed by pairs $(r,s)$ with $1\le r,s\le k$.
So think of the columns of $\mu$ as labeled either by a group element $g$
or a pair of indices $(r,s)$.

The parametrization on the right is similar so $\mu$ also has $|G|$ rows
labeled by group elements and $k^2$ rows labeled by pairs of indices
(i,j).

If $g,h \in G$, and $1\le i,j,r,s\le k$, entries in $\mu$ might
be described as $\mu[g,h]$, $\mu[g;(r,s)]$, $\mu[(i,j);g]$ or
$\mu[(i,j);(r,s)]$.  It will be useful to break up $\mu$ into submatrices
$M$, $R$, $C$, and $N$ as follows:
$$\eqalign {
M_{|G|\times |G|}&:=\bigl(\mu[g,h]\bigr)_{g,h}\cr
R_{|G|\times k^2}&:=\bigl(\mu[g;(r,s)]\bigr)_{g,(r,s)}\cr
C_{k^2 \times |G|}&:=\bigl(\mu[(i,j);g]\bigr)_{(i,j),g}\cr
N_{k^2\times k^2}&:=\bigl(\mu[(i,j);(r,s)]\bigr)_{(i,j),(r,s)}\cr
}$$
Think of $\mu$ assembled from these pieces as 
$$\mu = \pmatrix{M& R \cr C&N \cr }$$ 

It will also be convenient to talk about isolated rows from $R$ and
columns from $C$, so define
$$\eqalign{
&r_g = \hbox{the $g$th row of $R$}\cr
\hbox{and}\quad&c_g = \hbox{the $g$th column of $C$}\cr}$$

\sh{ Remaining pentagons for products with three $m$'s}

The pentagons $mmmg/m$, $mmgm/m$, $mgmm/m$, and $gmmm/m$ are information
about symmetries inside of the big associator $\mu$.  Each one of the
four pentagons can be split into four statements, about the $M$, $R$,
$C$, and $N$ parts of $\mu$.  The next 16 equations are to hold for all
$a,b$, and $g$ in $G$.

Pentagon $mmmg/m$: 
$$\eqalignno{
\mu[a,b]\alpha_2(b,g)&=\mu[g\inv a,b]\alpha_1(g\inv a,g)\beta_3(g)&(11M)\cr
r_a&=\beta_3(g)r_{g\inv a}\left(\gamma_3(g)\inv\tens\id_k\right)&(11R)\cr
\alpha_2(a,g)c_a&=\left(\gamma_3(g)\tens\gamma_3(g)\right)c_a&(11C)\cr
\left(\gamma_3(g)\tens\gamma_3(g)\right)N&
  =N\left(\gamma_3(g)\tens\id_k\right)&(11N)\cr}$$
Pentagon $mmgm/m$:  
$$\eqalignno{
\mu[a,b]\beta_2(g,a)\beta_3(g)&=\mu[a,bg\inv]&(12M)\cr
r_a\left(\gamma_2(g)\tens\gamma_3(g)\inv\right)&=\beta_2(g,a)r_a&(12R)\cr
c_{g\inv a}&=\beta_3(g)\left(\id_k\tens\gamma_2(g)\right)c_a&(12C)\cr
\left(\id_k\tens\gamma_2(g)\right)N&
  =N\left(\gamma_2(g)\tens\gamma_3(g)\inv\right)&(12N)\cr}$$
Pentagon $mgmm/m$: 
$$\eqalignno{
\mu[a,b]\beta_2(g,b)&=\mu\left[g\inv a,b\right]&(13M)\cr
r_{g\inv a}&=r_a\left(\id_k\tens\gamma_2(g)\right)&(13R)\cr
\left(\gamma_2(g)\inv\tens\gamma_1(g)\right)c_a&=\beta_2(g,a)c_a&(13C)\cr
\left(\gamma_2(g)\tens\gamma_1(g)\inv\right)N&
  =N\left(\id_k\tens\gamma2(g)\right)&(13N)\cr}$$
Pentagon $gmmm/m$:
$$\eqalignno{
\mu\left[a,g\inv b\right]&=\mu[a,b]\alpha_2(g,b)&(14M)\cr
r_a\left(\gamma_1(g)\tens\gamma_1(g)\right)&=\alpha_2(g,a)r_a&(14R)\cr
c_{g\inv a}&=\left(\gamma_1(g)\tens\id_k\right)c_a&(14C)\cr
N\left(\gamma_1(g)\tens\gamma_1(g)\right)&
  =\left(\gamma_1(g)\tens\id_k\right)N&(14N)\cr}$$

\proclaim {4.6\qua Lemma}
For all $g,h\in G$, $\mu [g,h]=\mu[\epsilon,\epsilon]\alpha_2(h,h)\inv$. 
\endproclaim
\prf
Take $a=g$ and $b=\epsilon$ in (13M) to get
$\mu[g,\epsilon]=\mu[\epsilon,\epsilon]$ for all $g\in G$.  By (14M),
$\mu[g,h\inv h]=\mu[g,\inv h]\alpha_2(h,h)$; the claim follows.
\endprf

\sh{4.7\qua Proof of order-control theorem 1.1 }
By lemma 4.6, the $1\times |G|$ row vector $m_a = \bigl(\mu(a,b)\bigr)_b$
is actually independent of $a$ so we will write $\vec m$ for this
common value.

Equation (11R) implies that $r_g =
\beta_3(g)r_\epsilon\left(\gamma_3(g)\inv\tens\id_k\right)$ for all
$g\in G$.

Thus the first $|G|$ rows of $\mu$ look like
$$\pmatrix{
\vec m & r_\epsilon\cr
\vec m & \beta_3(g_2)r_\epsilon\left(\gamma_3(g_2)\inv\tens\id_k\right)\cr
\vec m & \beta_3(g_3)r_\epsilon\left(\gamma_3(g_3)\inv\tens\id_k\right)\cr
\vdots&\vdots\cr
\vec m & \beta_3(g_N)r_\epsilon\left(\gamma_3(g_N)\inv\tens\id_k\right)\cr}_{|G|\times (|G|+k^2)} $$
Since $\mu$ is invertible these are linearly independent.  If we write
$r_\epsilon = \sum_i u_i\tens v_i$, it follows that the $|G|\times(k+1)$
matrix
$$\pmatrix{1 & \sum_i u_i\cr
1 & \sum_i \beta_3(g_1)u_i\gamma_3(g_1)\inv \cr
\vdots&\vdots\cr
1 & \sum_i \beta_3(g_N)u_i\gamma_3(g_N)\inv u_i \cr}$$
has linearly independent rows; hence $|G|\le k+1$.\endprf

\proclaim {4.8\qua Proposition}
If $|G|=k+1$ then $\beta_3(h)=\alpha_2(h,h)$ for all $h\in G$.
\endproclaim

\prf
Similar to lemma 4.6, equations (11M) and (12M) imply that
$$\mu[g,h]=\mu[\epsilon,\epsilon]\beta_3(g)\beta_3(h)\inv\beta_2(h,g)\inv$$
for all $g,h\in G$.  Equating this with the expression from lemma 4.4,
we find
$$\mu[\epsilon,\epsilon]\beta_3(g)\beta_3(h)\beta_2(h,g)\inv =
  \mu[\epsilon,\epsilon]\alpha_2(h,h)\inv$$
Implicit in the proof of the order control theorem, if $|G|=k+1$ then
$\mu[\epsilon,\epsilon]$ cannot be zero (this would violate linear
independence).  So we can cancel and rearrange to get
$$\alpha_2(h,hg\inv)=\beta_3(h)\beta_3(g)\inv$$
(we have also used that $\alpha_2\ideq\beta_2$ and $\alpha_2$ is
bimultiplicative.) Now, let $g=\epsilon$ in the above.
\endprf

\section{Further reduction of pentagons when $|G|=k+1$}

We proceed in the maximal-group case; throughout this section, $G$
is a group of order $k+1$.  The results of the previous section are
summarized for this special case in the following result.

\proclaim {5.1\qua Proposition}
In a near-group category with fusion rule $(G,k)$ where $|G|=k+1$,
there exists a choice of bases such that:

{\rm (i)}\qua  $\alpha,\alpha_1,\alpha_2,\alpha_3,\beta_1,\beta_2,\beta_3$
are all identically equal to 1.

{\rm (ii)}\qua $\gamma_1,\gamma_2\inv$ and $\gamma_3$ are conjugate
representations of $G$ in diagonal matrices.

{\rm (iii)}\qua $\mu[a,b]$ is constant (independent of both $a$ and $b$).

{\rm (iv)}\qua If we define $\chi_i(g)=\gamma_1(g)_{i,i}$, for $i=1\ldots k$,
then the $\chi_i$ are precisely the $k$ nontrivial characters of $G$.
\endproclaim

\prf
All of (i) was proven in the previous section.  For $\beta_2,\alpha_2,
\beta_3$ vanishing use 4.6 and 4.8.

For (ii), we have already established that $\gamma_1$, $\gamma_2\inv$,
and $\gamma_3$ are conjugate representations; it remains to show that
all three can be simultaneously be diagonalized.  Define $\theta$ on
$G\times G\times G$ by
$$\theta(g,h,k)=\gamma_1(g)\gamma_2(h)\gamma_3(k)$$
Equations (4) through (6) and $\alpha_2$ vanishing imply that $\theta$
is a homomorphism; we can therefore choose basis of $\hom(m,m\tens m)$
so that $\theta$ splits as linear representations.  Then of course
$\gamma_1(g)=\theta(g,\epsilon,\epsilon)$ is diagonal for all $g$;
similarly $\gamma_2$ and $\gamma_3$.

Claim (iii) follows from 4.6 and $\alpha_2\ideq1$

Claim (iv) follows from linear independence in the proof of theorem 1.1.
\endprf

\proclaim {5.2\qua Corollary}
With respect to the basis given in 5.1, there is a permutation $\pi$
of order 3 in $S_k$ and constants $\xi(j)$ such that
$$\lambda(\epsilon)=\left(\xi(j)\delta_{i,\pi(j)}\right)_{i,j}$$
\endproclaim

\prf
Follows from 4.5 and 5.1 (iv).
\endprf

\sh {5.3\qua Remark}
To make this $\pi$ an invariant of the category (so it doesn't depend
on an ordered basis of $\hom(m,mm)$), we more correctly think of $\pi$
as a permutation on the set of nontrivial characters of $G$.

\sh{Notation}
By the numbering indicated in 5.1 (iv) we can think of an index $i$
in $\left\{1,2,\ldots,k\right\}$ as corresponding to one of the $k$
nontrivial characters of $G$.  If we let $\chi_\epsilon$ stand for the
trivial character of $G$ then the set
$\left\{\epsilon,1,2,\ldots,k\right\}$ becomes a group under $*$ where
we define
$$r*s=i \iff \chi_r\chi_s=\chi_i$$
By the definition of the $\chi's$, there is an equivalent
description using entries in $\gamma_1$ matrices:  For $r,s \in
\left\{1,2,\ldots,k\right\}$,
$$\eqalign{
r*s=i&\iff \gamma_1(g)_r\gamma_1(g)_s=\gamma_1(g)_i\quad\forall g\in G, i\in \left\{1,2,\ldots,k\right\}\cr
r*s=\epsilon&\iff\gamma_1(g)_r\gamma_1(g)_s=1\qquad\forall g\in G\cr}$$
If $r$ is an index then $r\inv$ will denote the inverse with respect to
this $*$ operation.

Of course the group we obtain this way is abstractly isomorphic to $G$
itself and if we choose an isomorphism then we can view the permutation
$\pi$ from 5.2 as acting on the nonidentity elements of $G$.

Precisely the same characters appear in the representation $\gamma_3$
although their order of appearence on the diagonal is permuted by
the $\pi$ of lemma 5.2.  We introduce an analogous operation $\circ$
on $\left\{\epsilon,1,2,\ldots,k\right\}$ by 
$$\eqalign{
r\circ s=i&\iff \gamma_3(g)_r\gamma_3(g)_s=\gamma_3(g)_i\quad\forall
  g\in G, i\in \left\{1,2,\ldots,k\right\}\cr
r\circ s=\epsilon&\iff\gamma_3(g)_r\gamma_3(g)_s=1\qquad\forall g\in G\cr}$$
The inverse of $r$ with respect to $\circ$ will be denoted $r\ci$.
There is a straightforward translation between the $*$ and $\circ$
operations via $\pi$:  If $r\circ s \ne \epsilon$, then
$$\eqalign { r\circ s &= \pi\left[\pi\inv(r)*\pi\inv(s)\right]\cr
\hbox{and}\quad r\ci &= \pi\left[\pi\inv(r)\inv\right]\cr}$$
The $\circ$ operation is therefore not strictly necessary but it
eliminates the frequent awkward expressions on the righthand side
of the above equations, considerably neatening statements such as the
following:

\proclaim {5.4\qua Proposition {\rm ($R$ and $C$ structure)}}
For any $g$ in $G$,

{\rm (i)}\qua $r_g[(r,s)]=0\iff r*s=\epsilon$

{\rm (ii)}\qua $c_g[(i,j)]=0\iff i\circ j = \epsilon$

\endproclaim

\prf For (i), first suppose $r*s=\epsilon$.  Then by (14N),
$N[(i,j);(r,s)]$ must be zero for all $(i,j)$.  Linear independence of
the columns of $\mu$ implies that $r_g[(r,s)]\ne 0$ for some $g\in G$,
but (11R) says that for any $g,h\in G$, $r_h[(r,s)]$ and $r_g[(r,s)]$
are related via multiplication by a unit.  So in fact $r_g[(r,s)]\ne 0$
for all $g\in G$.

Conversely, suppose $r_g[(r,s)]\ne0$.  Then by (14R),
$$r_g[(r,s)]=\gamma_1(g)_r\gamma_1(g)_s r_g[(r,s)]$$
for all $g\in G$.  Since $r_g[(r,s)]$ can be cancelled, we see that
$\gamma_1(g)_r\gamma_1(g)_s=1$ for all $g\in G$, which means by definition
$r* s = \epsilon$.

Proof of (ii) is similar, using equations (11N) and (11C) instead of
(14N) and (14C).
\endprf

The previous result suggests notational shorthand: since for a given $r$,
$r_g[(r,s)]$ is nonzero only for $s=r\inv$ we will abbreviate to $r_g(r)$
to refer to this nonzero entry.  Similarly, $c_g(i)$ will be short for
$c_g[(i,i\ci)]$.

\proclaim {5.5\qua Proposition {\rm($N$ structure)}}
$N[(i,j);(r,s)]\ne0\iff r*s=i$ and $i\circ j=r$.
\endproclaim

\prf
First suppose $N[(i,j);(r,s)]\ne 0$.  Then for all $g\in G$, (11N) implies
$$\gamma_3(g)_i\gamma_3(g)_j N[(i,j);(r,s)]=\gamma_3(g)_r N[(i,j);(r,s)]$$
and we can cancel the $N$ entry to obtain
$\gamma_3(g)_i\gamma_3(g)_j=\gamma_3(g)_r$ for all $g$, which is the
definition of $i\circ j=r$.  Similarly, equation (14N) implies $r*s=i$.

Conversely, suppose $r*s=i$ and $i\circ j=r$.  By 5.4, $r_g[(r,s)]=0$
for all $g\in G$, so linear independence of the columns of $\mu$ implies
that $N[(p,q);(r,s)]$ is nonzero for at least one pair of indices $(p,q)$.
But the first half of this proposition implies that the only possible
nonzero entry in this column is $N[(i,j);(r,s)]$.\break
\endprf

\sh{5.6\qua The pentagon $mmmm/g$}

This pentagon gives the following equations for all $a,b,g\in G$,
for all indices $1\le i,j,r,s\le k$.  They have been reduced using the
simplifications available from 5.1(i).  We write $N[-;(r,s)]$ for the
entire $(r,s)$ column of $N$, and $N[(i,j);-]$ for the $(i,j)$ row.
$$\eqalignno{
\sum_{c\in G}\mu[a,c]\mu[c,b]+r_a\left(\lambda(g)\tens\id_k\right)c_b &
  = \delta_{a,b\inv g}&(15M)\cr
\sum_{c\in G}\mu[a,c]r_c[(r,s)]+r_a\left(\lambda(g)\tens\id_k\right)
  N[-;(r,s)]&=0&(15R)\cr
\sum_{c\in G}c_c[(i,j)]\mu[c,b] + N[(i,j);-]\left(\lambda(g)\tens\id_k\right)
  c_b&=0&(15C)\cr
\sum_{c\in G}c_c[(i,j)]r_c[(r,s)]+N[(i,j);-]\left(\lambda(g)\tens\id_k\right)
  N[-;(r,s)]&\phantom{xxxxxxxxxxxxxxx}\hfill&\cr
&\hskip-1in=\lambda(g)[i,s]\lambda(g)[j,r]&(15N)\cr}$$

\proclaim {5.7\qua Lemma} For any $g\in G$:

{\rm (i)}\qua $r_g[(r,s)]=\gamma_3(g)_r\inv r_\epsilon[(r,s)]$

{\rm (ii)}\qua $c_g[(i,j)]=\gamma_1(g)_i\inv c_\epsilon[(i,j)]$

\endproclaim

\prf
Immediate from (11R) and (14C).
\endprf

\proclaim {5.8\qua Proposition} The permutation $\pi$ from 5.2 satisfies
$\pi(x\inv)*\pi\inv(x)=\epsilon$ (for any index $x$, $1\le x \le k$.)
\endproclaim

\prf
In equation (15R) the summation in fact vanishes.  To see this,
fix $a,r,$ and $s$.  Now, $\mu[a,c]$ is constant by 5.1(iii), and
$r_c[(r,s)]=r_\epsilon[(r,s)]\gamma_3(c)_r\inv$ by 5.7.  Therefore,

$$\sum_{c\in G}\mu[a,c]r_c[(r,s)]=\mu[\epsilon,\epsilon]r_\epsilon[(r,s)]
  \sum_c\gamma_3(c)_r\inv=0$$
by orthogonality of characters.
So (15R) reduces to
$$r_a\left(\lambda(g)\tens\id_k\right)N[-;(r,s)]=0$$
and in particular, for $g=a=\epsilon$,
$$r_\epsilon\left(\lambda(\epsilon)\tens\id_k\right)N[-;(r,s)]=0$$
By 5.5 this yields no useful information if $r*s=\epsilon$ so assume
$r*s\ne\epsilon$.
By 5.5, there is a single nonzero entry in $N[-;(r,s)]$, namely
$N[(i,j);(r,s)]$ where $r*s=i$ and $i\circ j=r$.  The lefthand side will
not be zero if $\pi$ manages to ``line up'' that nonzero entry from $N$
with one of the nonzero entries in $r_\epsilon$.  But we know where the
nonzero entries in $r_\epsilon$ live by 5.4.  So if (15R) is to be
satisfied there must not exist any triple of indices $r,s,x$ with
$r*s\ne\epsilon$, $\pi(r*s)*x=\epsilon$, and $(r*s)\circ x = r$.

For any index $r$, there is an $s$ so that $x\inv=\pi(r*s)$.  So (15R)
says $\pi\inv(x\inv)\circ x \ne r$ for any index $r$; by elimination it
must be that $\pi\inv(x\inv)\circ x = \epsilon$.  Now,
$$\eqalign{ &\pi\inv(x\inv)\circ x = \epsilon\cr
  \implies &\pi^{-2}(x\inv)*\pi\inv(x)=\epsilon\cr
  \implies &\pi(x\inv)*\pi\inv(x)=\epsilon\quad
    \hbox{(since $\pi^3=\hbox{id}$)\ }\cr}$$
\vglue -20pt\qed

\proclaim {5.9\qua Proposition} The permutation $\pi$ satisfies
$\pi(s*t)=\pi(t)*\pi\bigl[\pi(s)\inv*\pi(t\inv)]$ (for all $s\ne t\inv$)
\endproclaim

\prf
Consider $(15N)$ with $r\ne s\inv$, $i=\pi(s)$, and $j=\pi(r)$, so that
the summation on the left vanishes but the righthand side is nonzero.
By the $N$ structure proposition, the only nonzero entry in column $(r,s)$
of $N$ occurs in row $(r*s,x)$ where $x=(r*s)\ci\circ r$.  Expressed
purely in terms of the $*$ operation,
$$x=\pi\left[\pi\inv(r*s)\inv*\pi\inv(r)\right]$$
But for the lefthand side of $(15N)$ to be nonzero, it must be that
$i=\pi(r*s)*x$; since we have chosen $i=\pi(s)$, the following relation
appears:
$$\pi(s)=\pi(r*s)*\pi\left[\pi\inv(r*s)\inv*\pi\inv(r)\right]\quad
  (\forall r\ne s\inv)$$
The identity claimed in the proposition follows from a change of variable:
set $t=(r*s)$ and rewrite.  The change of variable is reversible so the
two identities are in fact equivalent.
\endprf

\sh {5.10\qua Proof of intermediate theorem 4.2}
The results in 5.2, 5.8, and 5.9 establish Theorem 4.2: If a monoidal
structure exists, then we can extract a permutation $\pi$ satisfying
the three indicated properties.\endprf

In section 6 we will complete Theorem 1.2 by identifying which groups
support such a $\pi$ but the remainder of this section will continue
the analysis of the pentagon equations so that we can give an explicit
matrix solution to the pentagon equations in 5.18.

\sh{5.11\qua Identifying primitive data}
At this point we have reduced the problem of describing the categorical
data to describing the associativities $\lambda(\epsilon)$ and $\mu$.
We know the $\alpha$'s and $\beta$'s are trivial and we know that the
$\gamma$'s contain all the nontrivial characters.  All $\lambda(g)$
are determined once we know $\lambda(\epsilon)$ and the $\gamma$'s.
The $\lambda(\epsilon)$ is described by the permutation $\pi$ and the
function $\xi$ of 5.2; $\mu$ is described by $M$, $R$, $C$, and $N$ but
$M$ is a constant block; all of $C$ can be generated given the function
$c_\epsilon$.  Likewise $R$ is described by the function $r_e$, and $N$
by a function of two indices: $N(r,s)$, denoting the single nonzero
entry in column $(r,s)$ of $N$ (where $r*s\ne\epsilon)$.  Our next step
is to give a formula for $r_\epsilon$ in terms of $\ce$ thus reducing
the problem to specifying three functions $\xi$, $c_\epsilon$, and $N$.

\proclaim {5.12\qua Proposition}
$\re(i)=\left[|G|\xi(\pi\inv(i))\ce\left(\pi\inv(i)\right)\right]\inv$
\endproclaim

\prf
The key is to reformulate (15M) as a small matrix equation which
essentially allows us to solve for R in terms of C.  For $g\in G$, set
$r_g'=(r_g(i))_i$, a $1\times k$ row vector regarded as $\re$ with all
the zeros squeezed out. Similarly, set $c_g'=(c_g(i))_i^T$, a $k\times 1$
column vector.  Define
$$R_0 = \pmatrix{ 
y & \re' \cr
y & {r_1}' \cr
\vdots & \vdots \cr
y & {r_k}' \cr} \quad
C_0 = \pmatrix {
y & y & \cdots & y \cr
\ce' & c_1' & \cdots & c_k' \cr} \quad
\lambda_0 = \pmatrix {
1 & 0 \cr
0 & \lambda(\epsilon) \cr} $$
where $y=\sqrt{|G|\mu[\epsilon,\epsilon]^2}$.  The three matrices defined
above are $(k+1)\times(k+1)$ and invertible by linear independence of
characters.  For $g=\epsilon$ it is simple to verify that equation
(15M) for all $a,b \in G$ is equivalent to the matrix equation
$$R_0\lambda_0 C_0 = \left(\delta_{a,b\inv}\right)_{a,b\in G}$$
Write $X$ for the involutive permutation matrix
$\left(\delta_{a,b\inv}\right)_{a,b\in G}$ and $D_0$ for a diagonal matrix
with $(\ce(i)\ce(i\inv)|G|)\inv$ on the diagonal.  Character relations
let us calculate
$$C_0\inv = C_0^T X D_0$$
and so
$$R_0 = X C_0^T X D_0 \lambda_0\inv$$
Now it is easy to expand the matrix product on the right hand side and
obtain the formula claimed.
\endprf
Moreover, examination of the first column on both sides implies
$y^2|G|\inv = 1$, and we have

\proclaim {5.13\qua Corollary}
$\mu[\epsilon,\epsilon]=\pm1/|G|$
\endproclaim
We introduce $\delta:=|G|\mu[\epsilon,\epsilon]=\pm1$ as this will occur
frequently in subsequent equations.

\sh{5.14\qua The big pentagon: $mmmm/m$}
Graphically, the bases for the full left and full right associated
product have the form
$$\ltreeb m{}{}{}y{}xm\qua
  \raise 20pt \hbox{\rm or}\qua
  \rtreeb m{}{}{}{\hat y}{}{\hat x}m\qua
  \raise 20pt \hbox{\rm respectively.}$$
On the left, for example, there are $k|G|$ basis elements with $x\in
G$, $k|G|$ basis elements with $y\in G$, and $k^3$ with $x=y=m$, for a
total of $k^3+2k|G|$ elements.  Rather than view this pentagon as one
large $(k^3+2k|G|)$--dimensional matrix equation, we break it down into
several more sensible submatrix equations based on this grouping of
basis elements:
$$\bordermatrix{
     &{x \in G \atop y=m}&{x=m\atop y\in G}&x=y=m\cr
\vbox to 16pt{}\hfill {\hat x\in G \atop \hat y=m}& (16)&(18)&(23)\cr
\vbox to 20pt{}\hfill {\hat x=m \atop \hat y\in G}& (17)& \hbox{---}
  &(20),(21)\cr
\vbox to 20pt{} {\hat x=\hat y=m} & (22) & (18),(19) & (24)\cr }$$
The table indicates which submatrix of the pentagon equation corresponds
to which of the following constraints.  We have used all reductions
available to express these in the minimal number of variables: these are
essentially functional equations in the functions $\ce$, $\xi$, and $N$,
and they are to hold for all indices $i,j$ for which they make sense:
$$\eqalignno {
\ce(i)\ce(\pi\inv(i))\inv&=\delta
  \xi(i\inv)\xi(\pi(i\inv))\xi(\pi\inv(i))&(16)\cr
\xi(i)\ce(i)&=\delta\xi(\pi(i\inv))\ce(\pi(i\inv))&(17)\cr
\ce(i)&=\delta \ce(\pi(i)\inv)&(18)\cr
\ce(i)&=\ce(j)N(i*j\inv,j)\cr&\qquad N(\pi(j)*\pi(i)\inv,\pi(j)\inv)&(19)\cr
\xi(i)\ce(i)&=\delta\xi(\pi(i\inv))\ce(\pi(i\inv))&(20)\cr
\xi(\pi\inv(i)*\pi(j\inv))\ce(\pi\inv(i)*\pi(j\inv))
  &=\xi(\pi\inv(i))\ce(\pi\inv(i))\cr&\qquad N(i,i\inv*j)N(i\inv,j)&(21)\cr
\ce(i)N(i,j)&=\xi(j)\ce(i*j)N(\pi(i*j)\inv,\pi(j))  &(22)\cr
\xi(\pi\inv(i))N(i,j)\ce(\pi\inv(i))
	&=\xi\left(\pi(\pi\inv(j)*\pi(i)\right)\xi\left(\pi\inv(i*j)\right)
	\inv)&\cr
&\qquad N(j,(i*j)\ce(\pi\inv(i*j))&(23)\cr
\ce\left((i*j)\inv\right)\ce(j)\inv\xi(j)\inv&
  =N\left(\pi\inv(i)*\pi(j),\pi(j)\inv\right)\cr
  &\qquad N\left(\pi\inv(i),\pi(j)\right)&(24)\cr}$$
{\bf Note}\qua The $y,\hat y\in G$ piece of the pentagon is automatically
satisfied by the properties of the permutation $\pi$, without introducing
any new relations.  This will be illustrated in section 5.18 below.
In general the process of transcribing the above equations from the
abstract pentagon is lengthy but requires no sophistication.  We will
illustrate the derivation of two of the equations and omit the similar
details of the remainder.

\sh{5.15\qua Derivation of equation (16)}

As the table indicates, this equation is supposed to come from the
following piece of the pentagon:
$$\ltreeb mi{}{}m{}gm\qua
  \raise 30pt \hbox{$\labmapsto{}$}\qua \rtreeb mj{}{}m{}hm$$
For each $g,h \in G$ we get a $k\times k$ submatrix of coefficients indexed by $i$ and $j$.  Using the reductions
available it is not difficult to compute these coefficients.

Following the two-step path in the pentagon, the first step is 
$$\ltreeb mi{}{}m{}gm\qua
  \raise 30pt \hbox{$\labmapsto{}$}\qua
  \ctreec mijmmpm$$
$$\eqalign {\bas m g m \tens &\bas g m m \tens {\bas m m m}_i
\cr \ato \sum_p & c_g[p,j] {\bas m m m}_p \tens {\bas m m m}_i \tens
  {\bas m m m}_j +\hbox {irrelevant terms}\cr}$$
and the second step is
$$\ctreec mijmmpm\qua
  \raise 30pt \hbox{$\labmapsto{}$}\qua
  \rtreeb mj{}{}m{}hm$$
$$\eqalign {
{\bas m m m}_p &\tens {\bas m m m}_i \tens {\bas m m m}_j \cr 
& \ato r_h[p,i] \bas m m h \tens \bas h m m \tens {\bas m m m}_j +
  \hbox{irrelevant terms}\cr}$$
Therefore the $(j,i)$ coefficient is $\sum_p c_g[p,j]r_h[p,i]$.  By 5.4
this sum will be zero unless there is an index $p=j\ci=i\inv$
which happens iff $j=\pi\inv(i)$, in which case the coefficient is
$c_g(i\inv)r_h(i\inv)$, and so the submatrix here is
$$\biggl( c_g(i\inv)r_g(i\inv)\delta_{j,\pi\inv(i)} \biggr)_{j,i}$$
Now we follow the three-step path in the pentagon; the first step is
$$\ltreeb mi{}{}m{}gm\qua
  \raise 30pt \hbox{$\labmapsto{}$}\qua
  \ltreed mlmgm$$
$$\eqalign {\bas m g m \tens &\bas g m m \tens {\bas m m m}_i\cr
&\ato \sum_l\lambda(g)[l,i] \bas m g m \tens \bas g m m \tens
  {\bas m m m }_l\cr}$$
From here, 
$$\ltreed mlmgm\qua
  \raise 30pt \hbox{$\labmapsto{}$}\qua
  \rtreed mlmhm$$
$$\eqalign {\bas m g m &\tens \bas g m m \tens {\bas m m m }_l\cr
&\ato \mu[h,g] \bas m m h \tens \bas h m m {\bas m m m}_l \pit \cr}$$
Note that $\mu [h,g]=\mu[\epsilon,\epsilon]$, independent of $h$ and $g$.
In the third association,
$$\rtreed mlmhm\qua
  \raise 30pt \hbox{$\labmapsto{}$}\qua
  \rtreeb mj{}{}m{}hm$$
$$\eqalign {\bas m m h &\tens \bas h m m {\bas m m m}_l \cr
&\ato \lambda(h)[j,l] \bas m m h \tens \bas h m m \tens
  {\bas m m m}_j \pit\cr}$$
so the submatrix from this path is 
$$\eqalign{\mu[\epsilon,\epsilon]&\lambda(h)\lambda(g)\cr
  &=\mu[\epsilon,\epsilon]\lam\gamma_1(h)\lam\gamma_1(g)\quad
  \hbox{by 4.5}\cr}$$
and the requirement of the pentagon equation is that
$$\eqalignno{
\mu[\epsilon,\epsilon]\lam\gamma_1(h)\lam\gamma_1(g)&
  =\biggl( c_g(i\inv)r_h(i\inv)\delta_{j,\pi\inv(i)} \biggr)_{j,i}&(16')\cr}$$

\proclaim {Claim} $(16')$
for arbitrary $g$ and $h$ follows from the case $g=h=\epsilon$
\endproclaim
\prf
The matrix product on the left is 
$$\biggl(\mu[\epsilon,\epsilon]\chi_i(g)\xi(i)\chi_{\pi(i)}(h)\xi(\pi(i))
  \delta_{j,\pi\inv(j)}\biggr)_{j,i}$$
So $(16')$ can be written as 
$$\eqalignno{\mee\xi(i)\xi(\pi(i))\chi_i(g)\chi_{\pi(i)}(h) &
  = c_g(i\inv)r_h(i\inv) \quad\forall i&(16'')}$$
But $c_g(i\inv)=\ce(i\inv)\chi_i(g)$ by $(14C)$.  Similarly
$r_h(i\inv)=\re(i\inv) \chi_{\pi(i)}(h)$.  All terms involving $g$
and $h$ in fact cancel, so the general case reduces to
$$\mee\xi(i)\xi(\pi(i))=\ce(i\inv)\re(i\inv)$$
as claimed.
\endprf
Now 5.12 replaces the $\re$ with an equivalent expression in $\ce$, 
$$\mee\xi(i)\xi(\pi(i))=\ce(i\inv)\left(|G|\xi(\pi\inv(i\inv))
  \ce(\pi\inv(i\inv))\right)\inv$$
Rearrange:
$$\ce(i\inv)=\delta\xi(i)\xi(\pi(i))\xi(\pi\inv(i\inv))\ce(\pi\inv(i\inv))$$
and we have shown that this piece of the pentagon is precisely equation
(16) (after a final change of variable $i\mapsto i\inv$).

\sh {5.16\qua Derivation of equations (18) and (19)}
Graphically, we are interested in the following basis elements:
$$\ltreeb m{}{}{}qimm\qua
  \raise 30pt \hbox{$\labmapsto{\alpha}$}\qua
  \rtreeb mr{}qmpmm$$
Following the two-step path in the pentagon, first
$$\eqalign{
{\bas m m m}_i\otimes &\bas m g m \otimes \bas g m m \cr
&\hskip-10pt\ato \gamma_1(g)[r,i] {\bas m g m} \otimes {\bas g m m}
  \otimes {\bas m m m }_r \pit\cr}$$
and then
$$\eqalign{{\bas m g m} \otimes &{\bas g m m} \otimes {\bas g m m }_r \cr
\ato & \mu[(p,q);g] {\bas m m m}_p \otimes {\bas m m m}_q
  \otimes {\bas m m m }_r\pit\cr}$$
So the total coefficient from this path is simply
$$\mu[(p,q);g]\gamma_1(g)[r,i]$$
Following the three-step path, first
$$\eqalign{
{\bas m m m}_i\otimes &\bas m g m \otimes \bas g m m \cr
&\ato \sum_{h\in G}\mu[h;g] {\bas m m m}_i \otimes \bas m m h
  \otimes \bas h m m \cr
&\phantom{\ato}\quad+\sum_{j,l} \mu[(j,l);g] {\bas m m m}_i
  \otimes {\bas m m m}_j \otimes {\bas m m m}_l \cr}$$
then
$$\eqalign{
{\bas m m m}_i &\otimes \bas m m h \otimes \bas h m m \cr
&\ato \gamma_2(h)[p,i] {\bas m m m}_p \otimes {\bas m h m}
  \otimes {\bas h m m} \cr}$$
and
$$\eqalign{
  {\bas m m m}_p &\otimes {\bas m h m} \otimes {\bas h m m} \cr
&\hskip-20pt\ato \mu[(q,r);h] {\bas m m m}_p \otimes {\bas m m m}_q
  \otimes {\bas m m m }_r\pit\cr}$$
while
$$\eqalign{
{\bas m m m}_i &\otimes {\bas m m m}_j \otimes {\bas m m m}_l \cr
&\hskip-25pt\ato \sum_s \mu[(p,s);(i,j)] {\bas m m m}_p
  \otimes {\bas m m m}_s \otimes {\bas m m m}_l \hbox{+irrel. terms}\cr}$$
and
$$\eqalign{
{\bas m m m}_p &\otimes {\bas m m m}_s \otimes {\bas m m m}_l\cr
&\hskip-34pt\ato \mu[(q,r);(s,l)] {\bas m m m}_p \otimes {\bas m m m}_q
  \otimes {\bas m m m}_r \pit\cr}$$
So the total coefficient from this path is 
$$\eqalign{&\sum_{h\in G}\mu[(q,r);h]\gamma_2(h)[p,i]\mu[h;g]\cr
  &\quad +\sum_{j,l,s} \mu[(q,r);(s,l)]\mu[(p,s);(i,j)]\mu[(j,l);g]\cr}$$
and the requirement of the pentagon is that 
$$\eqalign {\mu[(p,q);g]\gamma_1(g)[r,i] = &\sum_{h\in G}\mu[(q,r);h]
  \gamma_2(h)[p,i]\mu[h;g]\cr
&\quad +\sum_{j,l,s} \mu[(q,r);(s,l)]\mu[(p,s);(i,j)]\mu[(j,l);g]\cr}$$
We will handle this in two cases.

{\bf Case I}\qua $i=p$\qua  
In this case, the second summation on the righthand side
vanishes, by the $N$ structure theorem (because $i*j\ne i=p$).  Again,
$\mu[g,h]=\mee$ independent of $g$ and $h$, so we are left with
$$c_g(i,q)\gamma_1(g)[r,i]=\mu[\epsilon,\epsilon]\sum_h
  c_h(q,r)\gamma_2(h)[i,i]$$
Now we can use (14C) and (12C) to replace the $c_g$ and $c_h$ terms
with expressions in $\ce$:
$$\ce(i,q)\gamma_1(g\inv)[i,i]\gamma_1(g)[r,i]=\mee\sum_h\ce(q,r)
  \gamma_2(h\inv)[r,r]\gamma_2(h)[i,i]$$
If $r\ne i$ the lefthand side vanishes since $\gamma_1$ is diagonalized,
and the righthand side vanishes by orthogonality of characters.  So we
can reduce to
$$\ce(i,q)=\mee|G|\ce(q,i)$$
and both vanish by the $C$ structure theorem unless $q=i\ci$, so the
only requirement here is
$$\ce(i)=\delta\ce(i\ci)\hbox{\ for all\ }i$$
which is precisely equation (18).

{\bf Case II}\qua $i\ne p$\qua  In this case the first summation vanishes because
$\gamma_2$ is diagonalized, and we are left with
$$\eqalignno{c_g(p,q)\gamma_1(g)[r,i]=\sum_{j,l,s}
  \mu[(q,r);(s,l)]\mu[(p,s);(i,j)]\mu[(j,l);g]&&(25)\cr}$$
Note that for many choices of indices $p,q,i,r$ the equation will be
satisfied without introducing any new constraints since both sides vanish.
We need to identify precisely when this happens.  Fix indices $p,q,i,r$.

By the $C$--structure theorem and diagonalization of $\gamma_1$, the
lefthand side is nonzero iff
$$\eqalignno {p=q\ci\hbox{\ and\ } r=i&&(26)\cr}$$
By the $C$-- and $N$--structure theorems, the righthand side is nonzero
iff there exist indices $s,l$ satisfying
$$\eqalignno {q\circ r &= s&{\rm(27i)}\cr
s*j\ci&=q&{\rm(27ii)}\cr
p\circ s&=i&{\rm(27iii)}\cr
i*j&=p&{\rm(27iv)}\cr}$$
{\bf Claim:}\qua   The conditions of (26) are equivalent to the conditions of (27).

\prf
Assume (26), guaranteeing the lefthand side is nonzero.  In order to
satisfy (27iii) and (27iv) we must take
$$j=p*i\inv\hbox{\ and \ }s=p\ci\circ i$$
This $s$ and $j$ clearly satisfy (27i); with a little more work we verify
that (27ii) also holds:
$$\eqalign{s*j\ci &= (p\ci \circ i) * (p*i\inv)\ci\cr
&=\pi(\pi(p\inv)*\pi\inv(i))*\pi(p*i\inv)\inv\cr
&=\pi(p)\inv\quad\hbox{using 5.9}\cr
&=p\ci = q\cr}$$
as required.  Conversely, supppose the conditions of (27) are satisfied.
Then we simply repeat the above calculation to show that $q=p\ci$.
It follows from (27i) and (27iii) that $i=r$, so both conditions of (26)
are satisified.
\endprf

Therefore in considering the requirements of $(25)$ we only need the
case $p=q\ci$ and $r=i$, reducing to
$$c_g(p)\gamma_1(g)[i,i]=N[p\ci\circ i,(p\ci\circ i)\inv * p\ci)]
  N[i,i\inv*p] c_g(p*i\inv)$$
Using (14C) and cancelling $\gamma$'s, this reduces further to
$$\ce(p)=N[p\ci\circ i,(p\ci\circ i)\inv * p\ci)]
  N[i,i\inv*p]\ce(p*i\inv)$$
which, after a change of variable, is precisely equation (19), and we
have completed showing that equations (18) and (19) completely describe
the requirements of this piece of the pentagon.\endprf

Derivation of the remaining equations in the list is very similar to the examples we have worked.
The $y,\hat y\in G$ equation is a little different; we have claimed it is automatically satisfied by the
properties of $\pi$ already established, and we will now prove that.

\sh {5.17\qua The $y,\hat y\in G$ piece of the pentagon }

Here we consider the following entries in the pentagon:
$$\ltreeb m{}{}{}qimm\qua
  \raise 30pt \hbox{$\labmapsto{\alpha}$}\qua
  \rtreeb m{}{}{}hjmm$$
It is easy to verify that the two-step path of the pentagon gives zero
coefficient for this entry.  Walking through the three-step path as we
did for previous pieces gives a total coefficient of
$$\eqalign{&\sum_{a\in G} \mee^2 \gamma_2(a)[i,j]\cr
&+\sum_{p,q,l}\mu[h;(l,q)]\mu[(j,l);(i,p)]\mu[(p,q);g]\cr}$$
Now the first sum (over $a\in G$) is always vanishing, so the requirement
of the pentagon is that
$$\sum_{p,q,l}r_h(l,q)\mu[(j,l);(i,p)]c_g(p,q)=0$$
for all $i$ and $j$.  If $i=j$ the $\mu[(j,l);(i,p)]$ is always zero by
the $N$--structure theorem and the sum vanishes.  If $i\ne j$, there
could be at most one nonzero term;  by the $N$-- and $C$--structure
theorems this would occur for 
$$\eqalign{
l&=j\ci\circ i\cr
p&=i\inv *j\cr
\hbox{and}\quad q&=p\ci\cr}$$
But then 
$$\eqalign{
l*q&=(j\ci\circ i)*(i\inv*j)\ci\cr
&=\pi\left[\pi(j\inv)*\pi\inv(i)\right]*\pi\left[i\inv*j\right]\inv\cr
&=\pi(j)\inv\quad\hbox{by 5.9}\cr
&\ne \epsilon\cr}$$
so by the $R$--structure theorem $r_h(l,q)=0$ and in fact every term of
the summation is zero without imposing any further conditions.\endprf

\sh {5.18\qua Constructing matrix solutions to the pentagon equations}

Suppose we are given near-group fusion rule $(G,k)$, $|G|=k+1$, and we
are given a permutation $\pi$ on $G\setminus\epsilon$ satisfying the
three conditions of Theorem 4.2.  Guided by the preceding analysis,
we give a set of associativity matrices and prove that they satisfy all
pentagon equations.

Using the notation established earlier, set
$$\alpha,\alpha_1,\alpha_2,\alpha_3,\beta_1,\beta_2\beta_3 \equiv 1$$
Let $\chi_1,\ldots,\chi_k$ be the nontrivial characters of $G$, and choose
an identification of $G$ with its dual group of characters so that $\pi$
can be regarded as permuting the characters.

For indices $i,j,r$ write $i*j=r$ provided $\chi_i\chi_j=\chi_r$;
write $i*j=\epsilon$ provided $\chi_i\chi_j=\chi_\epsilon$,
the trivial character.  Let $i\inv$ stand for the inverse of $i$
with respect to this operation.  Write $i\circ j = r$ provided
$\pi\inv(i)*\pi\inv(j)=\pi\inv(r)$, and let $i\ci$ denote the inverse
with respect to this operation.  As indicated before, $i\ci = \pi(i)\inv$.

Now we can make the following definitions:
$$\eqalign{
\lambda(\epsilon) &= \bigl (\delta_{i,\pi(j)}\bigr)_{i,j}\cr
\gamma_1(g) &= \pmatrix {\chi_1(g)&&& \cr &\chi_2(g)&& \cr &&\ddots &
  \cr&&&\chi_k(g)\cr}\cr
\gamma_2(g) &= \lambda(\epsilon)\inv \gamma_1(g\inv) \lambda(\epsilon)\cr
\gamma_3(g) &= \lambda(\epsilon)\inv \gamma_2(g\inv) \lambda(\epsilon)\cr
\lambda(g)  &= \gamma_3(g)\lambda(\epsilon) \quad\hbox{for arbitrary
  $g\in G$}\cr}$$

Note that with these definitions, all $\gamma_2(g)$ and $\gamma_3(g)$ are diagonal, and
$$\eqalign{&\gamma_2(g)[i,i]=\chi_{\pi(i)\inv}(g)\cr
  \hbox{and\quad}&\gamma_3(g)[i,i]=\chi_{\pi\inv(i)}(g)\cr}$$
At this point we have defined all associativities except $\mu$ and it
clear that pentagon equations (1)--(10) of this paper are satisfied,
as well as (1)--(11) of \cite {TY} together, these account for all
four-term products involving two $m$'s or fewer).  The only associativity
left to define is $\mu$.

For $g,h\in G$ and $i,j,r,s \in \{1,2,\ldots,k\}$, define
$$\eqalign{\mu[g,h]&=|G|\inv\cr
  \mu[(i,j);g]&=\chi_i(g\inv)\delta_{i,j\ci}\cr
  \mu[g;(r,s)]&=\left(|G|\chi_{\pi\inv(r)}(g)\right)\inv\delta_{r,s\inv}\cr
  \mu[(i,j);(r,s)]&=\delta_{r,i\circ j}\delta_{i,r*s}\cr}$$
Equations (11M--14M) are clearly satisfied by construction.

(14C) is clearly satisfied by construction. For (13C), check:  
$$\bigl(\left(\gamma_2(g)\inv\tens\gamma_1(g)\right)c_g\bigr)[(i,j)]
  =\chi_{\pi(i)}(g)\chi_j(g) c_g(i,j)$$
so if $j\ne i\ci = \pi(i)\inv$ then $c_g(i,j)=0$ and (13C) is satisfied;
if $j=i\ci$ then $\chi_{\pi(i)}(g)\chi_j(g)=1$ and (13C) is still
satisfied.  Verifying (12C) and (11C) is similar.

(11R) is clearly satisfied by construction.

Check that (12R) is satisfied: for $a,g \in G$,
$$\eqalign{&r_a\bigl(\gamma_2(g)\tens\gamma_3(g)\inv\bigr)=r_a\cr
  \iff&\re[r,s]\chi_{\pi(r)\inv}(g)\chi_{\pi\inv(s)\inv}(g)=\re[r,s]\quad
  \forall r,s\cr}$$
If $r\ne s\inv$ then $\re(r,s)=0$, and (12R) is satisfied.  If $r=s\inv$
then by property (ii) of $\pi$, $\pi(r)\inv * \pi\inv(s)\inv = \epsilon$,
so $\chi_{\pi(r)\inv}(g)\chi_{\pi\inv(s)\inv}(g)=1$, and (12R) is satisfied.

Check that (13R) is satisfied: This is trivial for $r\ne s\inv$, so
assume $r=s\inv$.  By property (ii) of $\pi$, then, $\pi(s) =
\pi\inv(r)\inv$, so for any $a,g \in G$,
$$\eqalign{\bigl(r_a(\id_k\tens\gamma_2(g)\bigr)[r,s] &
  = \bigl(|G|\chi_{\pi\inv(r)}(a)\mu[(\pi\inv(r),s);\epsilon]\bigr)\inv
  \chi_{\pi(s)\inv}(g)\cr
  &=\bigl(|G|\chi_{\pi\inv(r)}(a)\chi_{\pi(s)}(g)
  \mu[(\pi\inv(r),s);\epsilon]\bigr)\inv\cr &
  =\bigl(|G|\chi_{\pi\inv(r)}(g\inv a)\mu[(\pi\inv(r),s);\epsilon]
  \bigr)\inv\cr
  &=r_{g\inv a}[r,s]}$$
so (13R) is satisfied.  Checking (14R) is similar to (12R), so details are omitted.

(11N) is satisfied provided 
$$\chi_{\pi\inv(i)}(g)\chi_{\pi\inv(j)}(g)N[(i,j);(r,s)]=
  \chi_{\pi\inv(r)}(g)N[(i,j);(r,s)]$$
for all $i,j,r,s$.  If $i\circ j \ne r$ then both sides
are zero, and (11N) is satisfied; if $i\circ j = r$ then
$\chi_{\pi\inv(i)}\chi_{\pi\inv(j)}=\chi_{\pi\inv(r)}$ so (11N) is
satisfied in that case as well.

(12N) is satisfied provided 
$$\chi_{\pi\inv(i)}(g)\chi_{\pi\inv(j)}(g)N[(i,j);(r,s)]=
  \chi_{\pi\inv(r)}(g)N[(i,j);(r,s)]$$
for all $i,j,r,s$.  If $N[(i,j);(r,s)]=0$ this is satisfied.
If $N[(i,j);(r,s)]\ne0$ then $r*s=i$ and $i\circ j = r$ so
$j=r\circ(r*s)\ci$; we use the hypotheses on $\pi$ to massage this
expression:
$$\eqalign {j=r\circ(r*s)\ci&=\pi\bigl(\pi\left((r*s)\inv\right)
  *\pi\inv(r)\bigr)\cr &=\pi(s)*\pi(r*s)\inv\cr}$$
therefore
$$\eqalign{\pi(j)\inv &= \pi\bigl(\pi(s)*\pi(r*s)\inv\bigr)\inv\cr
&=\bigl(\pi\inv(s)*\pi(r)\bigr)\inv\quad\cr
&=\pi(r)\inv*\pi\inv(s)\inv\cr}$$
therefore in this case
$\chi_{\pi(j)\inv}=\chi_{\pi(r)\inv}\chi_{\pi\inv(s)\inv}$ and equation
(12N) is satisfied.

Checking (13N) is similar to (12N) and checking (14N) is similar to (11N) so details are omitted.

Check for (15M): by construction,
$$(\lambda(g)\tens\id_k\bigr)[(i,j);(r,s)]=
  \chi_r(g)\delta_{i,\pi(r)}\delta_{j,s}$$
and
$$c_b[(r,s)]=\chi_i(b\inv)\delta_{i,\pi(j)\inv}$$
So we can compute
$$\eqalign{
\left(\lambda(g)\tens\id_k\right)[(r,s)]&=\sum_{u,v}\chi_u(g)
  \delta_{i,\pi(u)}\delta_{j,v}\chi_u(b\inv)\delta_{u,\pi(v)\inv}\cr
&=\chi_{\pi\inv(r)}(gb\inv)\delta_{r,s\inv}\cr}$$
whence
$$\eqalign{
r_a\left(\lambda(g)\tens\id_k\right)c_b&=\sum_{r,s}
  \left(|G|\chi_\pi\inv(r)(a)\delta_{r,s\inv}\right)\inv
\left(\chi_\pi\inv(r)(gb\inv)\delta_{r,s\inv}\right)\cr
&=\cases {(|G|-1)/|G|,&$a=gb\inv$ \cr
-1/|G|,&$\ne gb\inv$ \cr}\cr}$$
and
$$\sum_{c\in G}\mu[a,c]\mu[c,b]+r_a\left[\lambda(g)\tens\id_k\right]c_b =
  \delta_{a,b\inv g}$$
so (15M) is satisfied.

Verifications of (15C) and (15R) are similar and easy.  We will check that
(15C) is satisfied, splitting into two cases.

{\bf Case I}\qua $i=j\ci$\qua In this case $N[(i,j);(r,s)]=0$ for all $r,s$ so
$$N[(i,j);-]\left(\lambda(g)\tens\id_k)\right)c_b=0$$
As for the summation,
$\sum_{c\in G}c_c[(i,j)\mu[c,b]=\sum_{c}\chi_i(c\inv)(1/|G|)=0$, so
(15C) is satisfied in this case.

{\bf Case II}\qua $i\ne j\ci$\qua  In this case every term in the summation is zero.  The other expression on the left also vanishes:
$$N[(i,j);-]\left(\lambda(g)\tens\id_k)\right)c_b=
  \sum_{r,s}N[(i,j);(r,s)]\chi_{\pi\inv(r)}
  \left(gb\inv\right)\delta_{r,s\inv}$$
so the $\delta$ kills off every summand where $N[(i,j);(r,s)]$ is nonzero.
That completes the check that (15C) is satisfied and (15R) is obviously
very similar.

Checking (15N): Recall this equation requires that
$$\sum_{c\in G}c_c[(i,j)]r_c[(r,s)]+N[(i,j);-]
  (\lambda(g)\tens\id_k)N[-;(r,s)]=\lambda(g)[i,s]\lambda(g)[j,r]$$
holds for all $g\in G$ and all indices $i,j,r,s$.  This is the trickiest
equation to verify; we will first do a general calculation of both
expressions on the lefthand side and then proceed to verify equality in
several cases.

First calculation: 
$$\eqalign{\sum_{c\in G}c_c[(i,j)]r_c[(r,s)]
&={1\over|G|}\sum_c \chi_i(c\inv)\chi_{\pi\inv(r)}(c\inv)
  \delta_{i,j\ci}\delta_{r,s\inv}\cr
&=\delta_{i,j\ci}\delta_{r,s\inv}\delta_{i,\pi\inv(r)\inv}\cr}$$
Second calculation: By construction
$N[(i,j);-]\left(\lambda(g)\tens\id_k\right)N[-;(r,s)]$
will certainly be zero if $i=j\ci$ or $r=s\inv$.  Let us calculate the
value of this expression assuming that $i\ne j\ci$ and $r\ne s\inv$.
For any indices $p$ and $q$,
$$\eqalign{&\left(\left(\lambda(g)\tens\id_k\right)
  \left(N[-;(r,s)]\right)\right)[(p,q)]\cr
&\qquad=\sum_{u,v}\left(\lambda(g)\tens\id_k\right)
  [(p,q);(u,v)]N[(u,v);(r,s)]\cr
&\qquad=\sum_{u,v}\chi_u(g)\delta_{p,\pi(u)}\delta_{q,v}
  \delta_{u,r*s}\delta_{v,(r*s)\ci\circ r}\cr
&\qquad=\chi_{r*s}(g)\delta_{q,(r*s)\ci\circ r}
  \delta_{p,\pi(r*s)}\cr}$$
Therefore, 
$$\eqalign{N[(i,j);-]&\left(\lambda(g)\tens\id_k\right)N[-;(r,s)]\cr
&=\sum_{p,q}\chi_{r*s}(g)\delta_{q,(r*s)\ci\circ r}
  \delta_{p,\pi(r*s)}\delta_{p,i\circ j}\delta_{q,(i\circ j)\inv*i}\cr}$$
Now it's clear that the only potentially nonzero term in this sum occurs
at the indices
$$\eqalign{&p=\pi(r*s)\cr \hbox{and}\quad&q=(r*s)\ci\circ r\cr}$$
That term will actually be nonzero iff $i\circ j = \pi(r*s)$ and $(i\circ
j)\inv*i=(r*s)\ci$.  We can now solve for the unique $i$ and $j$ that
will make this sum nonzero, for a particular $r$ and $s$:
$$\eqalign{\pi(r*s)\inv*i&=(r*s)\ci\circ r\cr
&=\pi\bigl(\pi\left(\left(r*s\right)\inv\right)*\pi\inv(r)\bigr)\cr
&=\pi(r*s)\inv*\pi\left(r\inv*r*s\right)\quad\hbox{using hypotheses on $\pi$}\cr
&=\pi(r*s)\inv*\pi(s)\cr}$$
hence a nonzero result occurs iff $i=\pi(s)$ and $j=\pi(r)$.  Thus the
final calculation is
$$\eqalign{&N[(i,j);-]\left(\lambda(g)\tens\id_k\right)N[-;(r,s)]\cr
&\qquad=\cases{0,&$i=j\ci\hbox{\ or\ }r=s\inv$ \cr
\chi_{r*s}(g)\delta_{i,\pi(s)}\delta_{j,\pi(r)},&\hbox{otherwise} \cr}\cr}$$
Now we will proceed to verify equality in (15N).  First of all, the
righthand side is easy to calculate; by construction,
$$\lambda(g)[i,s]\lambda(g)[j,r]=\chi_r(g)\chi_s(g)
  \delta_{i,\pi(s)}\delta_{j,\pi(r)}$$
From here we proceed in cases:

{\bf Case I}\qua When the righthand side is nonzero, ie, $i=\pi(s)$
and $j=\pi(r)$.

{\bf Case I, part (i)}\qua  Assume $r*s=\epsilon$.  It follows that $i\circ j
=\epsilon$ as well; and on the righthand side of the equation (15N)
we have specifically $\chi_r(g)\chi_s(g)=1$.

On the lefthand side $r*s=\epsilon$ implies the term involving $N$'s is
zero, and by the first calculation above, the summation is 1, so (15N)
is satisfied.

{\bf Case I, part (ii)}\qua  Assume $r*s\ne\epsilon$.  It follows that $i\circ
j \ne\epsilon$ as well.  On the lefthand side of equation (15N), all
terms in the summation are zero, whereas by second calculation the term
involving $N$'s gives precisely $\chi_{r*s}(g)$, so (15N) is satisfied.

{\bf Case II}\qua When the righthand side is zero, ie, $i\ne\pi(s)$ or
$j\ne\pi(r)$.  Note in this case the term involving $N$'s on the lefthand
side is zero, so we have to check that the summation vanishes.

{\bf Case II, part (i)}\qua Suppose $i\ne j\ci$ or $r\ne s\inv$.  In this case
every term of the summation is clearly zero.

{\bf Case II, part (ii)}\qua Suppose $i=j\ci$ and $r\ne s\inv$.  In this
case it must be that both $i\ne\pi(s)$ and $j\ne\pi(r)$.  For suppose
$i=\pi(s)$.  Then $\pi(s)=j\ci=\pi(j)\inv$ so $\pi(j)*\pi(s)=\epsilon$
and it follows that $\pi(j)=\pi\inv(r)$, or $j=\pi(r)$ contradictory
to the hypothesis we are assuming in Case II.  Similarly $j=\pi(r)$
forces $i=\pi(s)$, a contradiction.  So we suppose both $i\ne\pi(s)$
and $j\ne\pi(r)$.  Then $i=j\ci=\pi(j)\inv\ne\pi\inv(r)\inv$ so the
$\delta_{i,\pi\inv(r)\inv}$ term in our calculation of the summation
kills it, and equality holds in (15N).

This completes the verification of (15N) in all cases; equivalently,
the verification of the pentagon equations for the group summands of
a fourfold product of $m$'s.  The only remaining pentagon to check is
the $m$ summands of fourfold product of $m$'s, which we have reduced to
equations (16)--(24).  However, our construction makes the functions $\ce$,
$\xi$, and $N$ which we use to express those equations all identically 1,
giving a trivial solution to equations (16)--(24).  With this observation
we complete the verification that all pentagon equations in the category
are satisfied and we have constructed a coherent monoidal
structure.\endprf

The example data for $k=1,2,3$ in section 3 shows that this ``trivial''
solution need not be the only one.  Completing the classification for
those small values of $k$ ad hoc is easy but so far a general and useful
classification theorem hasn't appeared, so for the present we will be
content with the existence theorem.

\section{Groups that support a $\pi$ }

In this section we characterize the groups and permutations satisfying
the hypotheses of Theorem 4.2.  For reference: the hypotheses are that
$\pi$ is a permutation on the nonidentity elements of $G$, satisfying

{\rm (i)}\qua $\pi^3=\hbox{id}$

{\rm (ii)}\qua $\pi(x)\inv = \pi\inv(x\inv)$

{\rm (iii)}\qua $\pi(st)=\pi(t)\pi\bigl[\pi(s)\inv\pi(t\inv)]$ (for all
$s\ne t\inv$)

The result is:

\proclaim {6.1\qua Theorem}
A finite abelian group $G$ admits a permutation $\pi$ satisfying the
above hypotheses if and only if $G$ is the multiplicative group of a
field; ie, $G$ is cyclic of order $(p^\alpha-1)$.
\endproclaim

First we prove the following:

\proclaim {6.2\qua Lemma}
Suppose $G,\pi$ satisfy the hypotheses.  Then there
is an element $\omega\in G$ such that $s\pi(s)\pi\inv(s)=\omega$ for
all $s\ne\epsilon\in G$.  Moreover, $\omega^2=\epsilon$.
\endproclaim

\prf
This is established by iterating the product expansion given by property
(iii) of $\pi$.  Properties (i) and (ii) are also used to simplify.
Let $s\ne t\inv$ be any nonidentity elements in $G$.
$$\eqalign{\pi(st)&=\pi(t)\pi\bigl(\pi(s)\inv\pi(t\inv)\bigr)\cr
&=\pi(t) s\inv \pi\bigl(\pi(t)\pi\inv(s)\bigr)\cr
&=\pi(t) s\inv \pi\inv(t) \pi(s\inv t\inv)\cr}$$
Since $\pi(st)$ is symmetric in $s$ and $t$, we have
$$\pi(t)s\inv\pi\inv(t) = \pi(s) t\inv \pi\inv(s)$$
or, on rearranging,
$$s\pi(s)\pi\inv(s)=t\pi(t)\pi\inv(t)$$
Since this holds for arbitrary $s$ and $t$ we conclude that the product
$s\pi(s)\pi\inv(s)$ is independent of $s$ and we call the common value
$\omega$.  (The case where $G$ has only two nonidentity elements requires
a slightly different argument; in that case it is easy to verify that the
only permutation satisfying the hypotheses is the identity permutation,
and the present proposition is trivial).

For $\omega^2=\epsilon$, let $s\ne\epsilon\in G$.  By the preceding,
$$\eqalign{s\pi(s)\pi\inv(s)&=s\inv\pi(s\inv)\pi\inv(s\inv)\cr
&=\left(s\pi\inv(s)\pi(s)\right)\inv\cr}$$
hence $\omega=\omega\inv$.
\endprf
{\bf Note}\qua  It may be the case that $\omega=\epsilon$.

Now we can proceed with the main result.

\rk{Proof of 6.1} First suppose $G=\fmul$, the multiplicative group of
a field.  Then defining $\pi(x)=(1-x)\inv$ for $x\ne1$ gives a permutation
satisfying the hypotheses (easily verified).

For the converse, suppose $G,\pi$ satisfy the hypotheses.  We construct
a field $F$ with multiplicative group $G$ as follows:  setwise, of course,
$F=G\cup\{0\}$.  The multiplication extends to $F$ by the obvious rule
$0*x=x*0=0$.  The main chore is to define addition.  We do this as
follows:
\nobreak
$$\eqalign{0+a=a+0&=a\cr
a+b&=\pi(\omega ba\inv)\inv a,\quad a,b\ne 0,b\ne \omega a\cr
a+b&=0,\quad b=\omega a\cr}$$
The element $\omega$, apparently, acts as a unary $(-)$ operator,
giving additive inverses.  We now have to show that $F$ is an abelian
group under this addition and that multiplication distributes.  We will
repeatedly use the lemma to replace expressions of the form $\pi(x)\inv$
with $\omega x \pi\inv(x)$.

(1)\qua  Addition is commutative:
$$\eqalign{a+b &= \pi(\omega ba\inv)\inv a\cr
&=\omega\ \omega ba\inv\ \pi\inv(\omega ba\inv) a\cr
&=\pi(\omega ab\inv)\inv b\cr
&=b+a\cr}$$
(2)\qua  Addition is associative:  If any of $a,b,c$ is 0 then
$(a+b)+c=a+(b+c)$ is trivial, so we assume $a,b,c\in G$.  Expanding the
left-associated expression,
$$\eqalign{(a+b)+c&=\pi\left(\omega c\left(a+b\right)\inv\right)\inv
  \left(a+b\right)\cr
&=\pi\left(\omega c \left( \pi(\omega ba\inv)\inv a\right)\inv\right)\inv
  \pi(\omega b a\inv)\inv a\cr
&=\pi\left(\omega c\ \pi\left(\omega b a \inv\right)a\inv\right)\inv\ \pi
  (\omega ba\inv)\inv a\cr}$$
On the other (the right) hand,
$$\eqalign{a+(b+c)&=\pi\left(\omega\left(b+c\right)a\inv\right)\inv a\cr
&=\pi(\omega\ \pi(\omega cb\inv)\inv\ ba\inv)\inv a\cr
&=\pi(st)\inv a\cr}$$
where $s=\pi(\omega cb\inv)\inv$ and $t=\omega ba\inv$.  Continuing the
calculation, we have
$$\eqalign{a+(b+c)&=\left(\pi(t)\pi
  \left(\pi(s)\inv\pi(t\inv)\right)\right)\inv a\cr
&=\pi\left(\pi(s)\inv\pi(t\inv)\right)\inv \pi(\omega ba\inv)\inv a\cr}$$
now we evaluate the expression $\pi(s)\inv\pi(t\inv)$:
$$\eqalign{\pi(s)\inv\pi(t\inv)&=\pi\bigl(\pi(\omega cb\inv)\inv\bigr)\inv
  \pi(\omega ab\inv)\cr
&=\pi\bigl(\pi\inv(\omega bc\inv)\bigr)\inv\pi(\omega ab\inv)\cr
&=\omega cb\inv \pi(\omega ab\inv)\cr
&=\omega cb\inv\ \omega\ \omega ba\inv\ \pi(\omega ba\inv)\cr
&=\omega c\pi(\omega ba\inv) a\inv\cr}$$
and putting it all together we have obtained
$$a+(b+c)=\pi\left(\omega c\ \pi\left(\omega b a \inv\right)
  a\inv\right)\inv\ \pi(\omega ba\inv)\inv a$$
which agrees with the expression obtained from left association.

At this point we have established that $F$ is an abelian group under $+$,
and we need only prove that multiplication distributes.  This is easy:
if $a,b,c\ne 0$,
$$\eqalign{a(b+c)&=a\ \pi(\omega cb\inv)\inv\ b\cr
&=\pi(\omega ac(ab)\inv)\inv ab\cr
&=ab+ac\cr}$$
while if any of $a,b,c$ is 0 the calculation is trivial.  Distributivity
on the right follows by commutativity of multiplication.  So we have $F$
a field with units the group $G$, as claimed.\endprf

Note that in the field $F$ we have constructed from $G$ and $\pi$, the
permutation $\pi$ can now be expressed by the formula $\pi(x)=(1-x)\inv$
so this is essentially the only example.

\section {The trivial group case: a nonexistence result}

Here we study the case in which our category has only two
isomorphism classes of simple objects: the identity $\epsilon$
and a noninvertible $m$; the fusion rule is described by $m\tens
m = \epsilon \oplus km$.  There are two nontrivial associativity
morphisms, $\lambda\co\hom(\epsilon,(mm)m)\to\hom(\epsilon,m(mm))$ and
$\mu\co\hom(m,(mm)m)\to\hom(m,m(mm))$.  These have to satisfy two pentagon
equations: $mmmm/\epsilon$ and $mmmm/m$.  Element-level descriptions of
these equations such as we produced in the nontrivial group case are not
particularly helpful since we can't use character relations to deduce
large amounts of symmetry.  We can, however, extract some information
by expressing the pentagons in the large as matrix equations and using
determinants.

The $mmmm/\epsilon$ pentagon is formulated as a $(k^2+1)\times (k^2+1)$
matrix equation.  Write $\id_k$ for a $k\times k$ identity matrix, and
$X_k$ for the matrix which operates on basis $\{e_i\tens e_j\}_{i,j=1\ldots k}$
by tensor flip: $X_k\co e_i \tens e_j \mapsto e_j\tens e_i$.  Write
$\lambda'$ for the block matrix $(1)\oplus (\lambda\tens\id_k)$.
With this notation, the pentagon equation is
$$\mu\lambda'\mu = \lambda' \bigl( (1)\oplus X_k \bigr) \lambda'$$
Write $M=\det \mu$ and $L=\det \lambda$.  Passing to determinants from
the above equation,
$$\eqalignno{M^2L^k&=L^{2k}\det(X_k)&(25)\cr}$$
The $mmmm/m$ pentagon is formulated as a $k^3+2k$ dimensional matrix
equation.  We need some permutation matrices to take care of reordering
bases in between steps of the pentagon (as the $X_k$ did in the smaller
pentagon).  Define
$$P_1 = \pmatrix {0 & \id_k & 0 \cr
\id_k & 0 & 0 \cr
0 & 0 & \id_{k^3} \cr}$$
in other words, $P_1$ exchanges the first block of $k$ basis elements
with the second $k$.  Clearly $\det P_1 = (-1)^k$ and $P_1^2=\id$. Define
$$P_2 = \id_{2k}\oplus
  \overbrace{X_k \oplus \cdots \oplus X_k}^{\hbox{$k$ copies}}$$
and note $P_2$ is also involutory, with $\det P_2 = (\det X_k)^k$.

Finally, let $P_3$ be the permutation which changes the basis elements 
$$\ltreeb mr{}qypxm$$
\noindent for $\hom(m,(((mm)m)m))$ from $(x,y,p,q,r)$ order to
$(x,p,y,q,r)$ order. Now, if we write
$$\eqalign{A&=\lambda \oplus\left(\id_k\tens \mu\right)\quad\hbox{and}\cr
  B&=\id_k\oplus\left(\mu\tens \id_k\right)\cr}$$
then the large pentagon $mmmm/m$ is expressed as the matrix equation
$$A^{P_3}B^{P_1}A^{P_3} = P_1BP_2P_1BP_1$$
or
$$A^{P_3}B^{P_1}A^{P_3} = B^{P_1}P_1P_2B^{P_1}$$
where exponents denote conjugation.

Now $\det A = L^kM$ and $\det B=L^k$ so on passing to determinants from
the pentagon, we obtain
$$\eqalignno{(LM^k)^2&=M^k(-1)^k(\det X_k)^k&(26)\cr}$$
\medskip
\proclaim {7.1\qua Lemma} 
$\det X_k=
\cases {+1,& $k\equiv 0\hbox{\rm\ or }1\pmod 4 $\cr
       -1,& $k\equiv 2\hbox{\rm\ or }3\pmod 4 $\cr}$
\endproclaim

\prf
The $(-1)$--eigenspace is spanned by $k(k-1)/2$ vectors of the form $e_i\tens e_j - e_j\tens e_i$, $1\le i < j \le k$.
The dimension will be even, hence $\det X_k=+1$, iff $4|k$ or $4|(k-1)$.
\endprf

\proclaim {7.2\qua Corollary}
If $k\equiv 2\hbox{\rm\ or }3\pmod 4$ then $L$ is a root of $-1$.
\endproclaim

\prf
{\bf Case I}\qua $k\equiv 2\pmod 4$\qua  
Write $k=2r$, where $r$ odd.  Equations 25
and 26 reduce to
$$\eqalign{
L^2&=M^{-2r}\cr
M^2&=(-1)L^{2r}\cr}$$
whence $L^{2r^2+2}=-1$.

{\bf Case II}\qua $k\equiv 3\pmod 4$\qua  In this case equations 25 and 26 reduce to
$$\eqalign{
L^2&=M^{-k}\cr
M^2&=(-1)L^k\cr}$$
whence $L^{k^2+4}=-1$.
\endprf

\proclaim {7.3\qua Lemma}
$\mee$ is invertible.
\endproclaim

\prf
This a consequence of duality assumptions (explained in the introduction).
We have assumed the existence of a map $\lambda_m$ such that
$$m\to\epsilon m\to(mm)m\smash{\mathop{\rightarrow}\limits^{\mu}}
  m(mm)\smash{\mathop{\hbox to 25pt{\rightarrowfill}}
  \limits^{1\tens\lambda_m}}m\epsilon\to m$$
is equal to the identity on $m$, which means $\lambda_m\mee=1$.
\endprf

\proclaim {7.4\qua Proposition}
  $L^3=1$ (independent of $k$).
\endproclaim
\prf
To prove this we will examine more carefully certain submatrices in the
big pentagon.

First, consider the $k\times k$ submatrix corresponding to:
$$\ltreeb mj{}{}m{}em\qua
  \raise 30pt \hbox{$\labmapsto{}$}\qua
  \rtreeb mi{}{}m{}em$$
Working through both sides of the pentagon, we find this piece gives
the equation
$$\sum_p\lambda[p,j]\mee\lambda[i,p]=\sum_p
  \mu[(p,i);\epsilon]\mu[\epsilon;(p,j)]$$
If we define $k\times k$ matrices $\mu_R$ and $\mu_C$ as follows:
$$\eqalign{
\mu_R[i,j]&=\mu[\epsilon;(i,j)]\cr
\mu_C[i,j]&=\mu[(j,i);\epsilon]\cr}$$
then the above equation can be reformulated as a matrix equation:
$$\eqalignno{\mu_C\mu_R &= \mee \lambda^2&(27)\cr}$$
Now we consider the $k\times k$ submatrix corresponding to:
$$\ltreeb mj{}{}m{}em\qua
  \raise 30pt \hbox{$\labmapsto{}$}\qua
  \rtreeb m{}{}{}eimm$$
This piece gives us
$$\sum_{p,q}\lambda[p,k]\mu[(i,q);\epsilon]\mu[\epsilon;(q,p)]=\mee\id_k$$
which we can reformulate as a matrix equation:
$$\eqalignno{\mu_C^T\mu_R\lambda &= \mee\id_k&(28)}$$
By previous proposition, $\mee$ is invertible.  So we can combine (27)
and (28), pass to determinants, and cancel to obtain
$$\det \lambda^2 = (\det \mu_C) (\det \mu_R) = \det\lambda\inv$$
or $L^3=1$ as claimed.
\endprf

\proclaim {7.5\qua Proof of theorem 1.4} \endproclaim
By 7.2, 7.4, and the hypothesis that the ground ring $R$ has
characteristic other than 2, the pentagon equations admit no solution
when $k\equiv 2\hbox{\rm\ or }3\pmod 4$.\endprf

\section {Additional structures I:  commutativity}

With sections 3 through 5 we have given quite a few examples, and a
general construction, of monoidal structures in the maximal-group case.
Topological applications (eg, link invariants) generally require at
least a commutative structure in addition to the monoidal structures,
and that means examining the hexagon equations for coherence of the
two structures.  In this section we will list the hexagon equations
that need to be satisfied, reduce them as much as possible in general,
and then show that the general construction of 5.18 always admits a
(symmetric) commuting structure.

Throughout this section we will assume we are in the setting of theorem
1.2 and have chosen bases so that all the reductions of sections 4 and
5 are in effect.

\sh{8.1\qua Notation for commutativity}

This is a slight extension of the notation used previously in \cite{S}.  
For $g,h\in G$, the commuting isomorphisms are:
$$\halign{\indent#&#&#\hfil\cr
$gh \to hg$&\quad is multiplication by &$\sigma_0(g,h)$ \cr
$gm \to mg$&\quad is multiplication by &$\sigma_1(g)$ \cr
$mg \to gm$&\quad is multiplication by &$\sigma_2(g)$ \cr
$mm \to mm$&\quad is multiplication by &$\sigma_3(g)$ on the $g$ summand\cr}$$
and commutativity on the $m$ summands of $m\tens m$ is represented by
a $k\times k$ matrix $\sigma_4$.

\sh {8.2\qua Unreduced hexagon equations }

The point here is simply to write down a transcription of the hexagon
equations using the notation we have developed, without attempting
any simplifications.  Once we have an entire transcription in hand,
we will proceed to analyze and simplify.

We will refer to the $xyz/w$ hexagon meaning the content of the hexagon
identity for the $w$ summands in the product $xyz$.  Note that lots of
associativities will be invisible since we have chosen bases well.

Let $a,b,c$ stand for arbitrary elements of $G$.

Hexagons for products with no $m$'s:
$$\eqalignno{\sigma_0(a,c)\sigma_0(a,b)&=\sigma_0(a,bc)&(abc/abc)\cr}$$
Hexagons for products with one $m$:
$$\eqalignno{\sigma_1(a)\sigma_0(a,b)&=\sigma_1(a)&(abm/m)\cr
  \sigma_0(a,b)\sigma_1(a)&=\sigma_1(a)&(amb/m)\cr
  \sigma_2(b)\sigma_2(a)&=\sigma_2(ab)&(mab/m)\cr}$$
Hexagons for the group summands of products with two $m$'s:
$$\eqalignno{\sigma_2(a)\sigma_3(ba\inv)&=\sigma_3(b)&(mma/b)\cr
  \sigma_3(ba\inv)\sigma_2(a)&=\sigma_3(b)&(mam/b)\cr
  \sigma_1(a)^2&=\sigma(a,ba\inv)&(amm/b)\cr}$$
Hexagons for the $m$ summands of products with two $m$'s:
\nobreak
$$\eqalignno{\sigma_2(a)\gamma_3(a)\sigma_4&=\gamma_2(a)
  \sigma_4\gamma_3(a)&(mma/m)\cr
  \sigma_4\gamma_1(a)\sigma_2(a)&=\gamma_1(a)\sigma_4\gamma_2(a)&(mam/m)\cr
  \sigma_1(a)\gamma_2(a)\sigma_1(a)&=\gamma_3(a)
  \sigma_1(a)\gamma_1(a)&(amm/m)\cr}$$
The hexagon for group summands in a product of three $m$'s:
$$\eqalignno{\sigma_4\lambda(g)\sigma_4&=\lambda(g)
  \sigma_3(g)\lambda(g)&(mmm/g)\cr}$$
The hexagon for $m$ summands in a product of three $m$'s we break into
individual entries in the matrices.

I\qua For all $g,h\in G$,
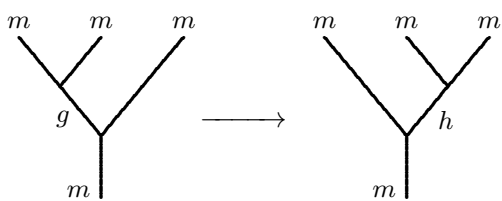
$$\ltreea mmm{}g{}m
  \raise 20pt \hbox{$\labmapsto{}$}\qua
  \rtreea mmm{}h{}m$$
$$\eqalign{\sigma_3(g)\mu[h,g]\sigma_3(h)&=
\sum_{a\in G}\mu[a,g]\sigma_2(a)\mu[g,a]\cr
&\quad+\sum_{i,j,p}\mu[(i,j);g]\sigma_4[p,i]\mu[h;(p,j)]\cr}$$
\bigskip
II\qua  For all $g\in G$, $i,j=1\ldots k$,
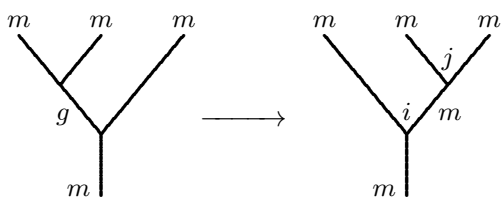
$$\ltreea mmm{}g{}m
  \raise 20pt \hbox{$\labmapsto{}$}\qua
  \rtreea mmmjmim$$
$$\eqalign{\sum_p\sigma_3(g)\mu[(i,p);g]\sigma_4[j,p]
&=\sum_{a\in G}\mu[a,g]\sigma_2(a)\mu[(i,j);a]\cr
&\quad+\sum_{r,s,t}\mu[(r,s);g]\sigma_4[t,r]\mu[(i,j);(t,s)]\cr}$$
\bigskip
III\qua  For all $g\in G$, $r,s=1\ldots k$,
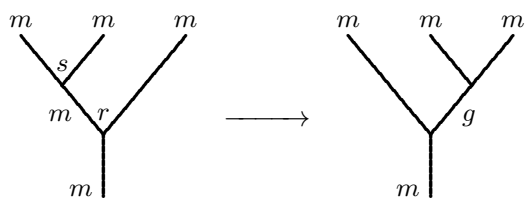
$$\ltreea mmmsmrm\qua
  \raise 20pt \hbox{$\labmapsto{}$}\qua
  \rtreea mmm{}g{}m$$
$$\eqalign{\sum_t \sigma_4[t,s]\mu[g;(r,t)]\sigma_3(g)
&=\sum_{a\in G}\mu[a;(r,s)]\sigma_2(a)\mu[g,a]\cr
&\quad+\sum_{j,p,q}\mu[(p,q);(r,s)]\sigma_4[j,p]\mu[g;(j,q)]\cr}$$
IV\qua  For all $i,j,r,s=1\ldots k$,
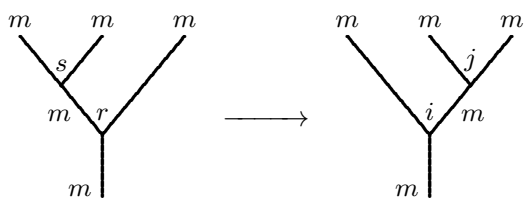
$$\ltreea mmmsmrm\qua
  \raise 20pt \hbox{$\labmapsto{}$}\qua
  \rtreea mmmjmim$$
$$\eqalign{\sum_{u,v}\sigma_4[u,s]\mu[(i,v);(r,u)]\sigma_4[j,v]&=
\sum_{a\in G}\mu[a;(r,s)]\sigma_2(a)\mu[(i,j);a]\cr
&\quad+\sum_{p,q,t}\mu[(p,q);(r,s)]\sigma_4[t,p]\mu[(i,j);(t,q)]\cr}$$
\bigskip
That completes the list of hexagon equations.  There are also inverse
hexagons to be considered since we do not assume a symmetric commutative
structure but we will postpone analysis of the inverse hexagons until
we have simplified the standard ones.

\sh{8.3\qua Reducing the hexagons }

Quite a lot of reductions are possible in the list of hexagons.  To start,
$\sigma_0\equiv 1$ is obvious.  We will show that $\sigma_1,\sigma_2$ must
both be equal to a distinguished linear character of $G$.  Then $\sigma_3$
is entirely determined up to the value of $\sigma_3(\epsilon)$, and
happily, the $\sigma_4$ matrix admits considerable simplification as well.

\proclaim{8.4\qua Proposition}
$\sigma_1=\sigma_2$
\endproclaim

\prf
It's clear from the list of hexagons that $\sigma_0\equiv 1$ so the
present proposition would follow from the $amm/b$ hexagon if we assumed
symmetric commutativity.  In general, however, it still follows from
the $mma/b$ hexagon together with its inverse:
$$\sigma_2(a)=\sigma_3(b)\sigma_3(ba\inv)\inv=\sigma_1(a)$$
From the $mab/m$ hexagon it's now clear that $\sigma_1$ and $\sigma_2$
are linear characters of $G$.  Our next result sharpens that statement.
By the results of section 7, we know there is a unique character
$\chi_\omega$ with $\chi_i\chi_{\pi(i)}\chi_{\pi\inv(i)}=\chi_\omega$
for all $i$; it satisfies $\chi_\omega^2=1$ and is the trivial character
(ie, $\omega=\epsilon$) iff $|G|$ is odd.
\endprf
 
\proclaim{8.5\qua Proposition} $\sigma_1=\sigma_2=\chi_\omega$ \endproclaim

\prf
By the $amm/b$ hexagon, $\sigma_1^2=1$.  Therefore the $amm/m$ hexagon
can be reduced to
$$\sigma_1(a)=\gamma_3(a)\gamma_1(a)\gamma_2(a)\inv$$
Since $\gamma$'s are diagonalized this is equivalent to
$$\sigma_1(a)=\chi_{\pi\inv(i)}(a)\chi_i(a)\chi_{\pi(i)}(a)$$
for all $a$ and $i$.  Hence, $\sigma_1=\chi_\omega$ as claimed.
\endprf

\proclaim{8.6\qua Proposition {\rm ($\sigma_4$ Structure)}}\hskip-5pt
There are invertible constants $\psi(1),\ldots,\psi(k)$ such
that the matrix $\sigma_4$ has the form $\sigma_4=\pmatrix{
\psi(j)\delta_{i,\pi(j\inv)} }_{i,j}$
\endproclaim

\prf
Consider the $mma/m$ hexagon.  It now reduces to
$$\chi_\omega(g)\sigma_4[i,j]=\chi_{\pi(i\inv)}(g)
  \chi_{\pi\inv(i\inv)}(g)\chi_{\pi\inv(j)}(g)\sigma_4[i,j]$$
for all $i,j,$ and $g$.  But that, together with invertibility of
$\sigma_4$ tells us that $\sigma_4[i,j]$ is nonzero (moreover, invertible)
iff $i=\pi(j\inv)$.
\endprf

Note that the $mam/m$ hexagon is satisfied by precisely the same condition.

\proclaim{8.7\qua Proposition}
All of the $mmm/g$ hexagons are satisified iff the constants $\psi$
satisfy the identity
$$\psi(j)\psi(\pi(j)\inv)\xi(\pi(j\inv))=\sigma_3(\epsilon)\xi(j)\xi(\pi(j))$$
for all $j$.
\endproclaim

\prf
This is simply working out the matrix products on both sides.  On the left,
$$\sigma_4\lambda(g)\sigma_4 = \biggl( \psi(j)\psi(\pi(j)\inv)
  \xi(\pi(j\inv))\chi_{\pi(j\inv}(g)\delta_{i,\pi\inv(j)} \biggr)_{i,j}$$
while on the right,
$$\lambda(g)\sigma_3(g)\lambda(g)=
\biggl( \xi(j)\xi(\pi(j))\sigma_3(\epsilon)\chi_\omega(g)
  \chi_j(g)\chi_{\pi(j)}(g)\delta_{i,\pi\inv(j)}\biggr)_{i,j}$$
All the $\chi$'s cancel and the proposition follows.
\endprf

\sh{8.8\qua Reducing the $mmm/m$ hexagon, part I}

On the lefthand side,
$$\sigma_3(g)\mu[h,g]\sigma_3(h)=\sigma_3(\epsilon)^2\chi_\omega(gh)\mee$$
On the righthand side, the first summation is
$$\sum_{a\in G}\mu[a,g]\sigma_2(a)\mu[g,a]=
  \left(|G|\inv\right)\sum_{a\in G}\chi_\omega(a)=
  |G|\inv\delta_{\omega,\epsilon}$$
The second summation is
$$\sum_{i,j,p}\mu[(i,j);g]\sigma_4[p,i]\mu[h;(p,j)]=
  \cases {|G|\inv\psi(\omega)\xi(\omega)\inv
  \chi_\omega(gh)\inv,&$\omega\ne\epsilon$ \cr
  0,&$\omega=\epsilon$\cr}$$
In case $|G|$ is odd, $\omega=\epsilon$, the requirement of this piece
of the hexagon is simply
$$\sigma_3(\epsilon)^2=\delta$$
(recall $\delta=|G|\mee=\pm1$).

In case $|G|$ is even, $\omega\ne\epsilon$, the requirement of the hexagon is
$$\sigma_3(\epsilon)^2=\psi(\omega)\xi(\omega)\inv$$

In either case this completes the simplification of this piece of the
$mmm/m$ hexagon.

\sh {8.9\qua Reducing the $mmm/m$ hexagon, part II}

On the lefthand side,
$$\eqalign{\sum_p \sigma_3(g)\mu[(i,p);g]\sigma_4[j,p]&=
  \sigma_3(g)\sum_p\ce[(i,p)]\chi_i(g)\inv\sigma_4[j,p]\cr
&=\sigma_3(\epsilon)\chi_\omega(g)\chi_i(g)\inv\ce(i)
  \psi\left(\pi(i)\inv\right)\delta_{i,\pi(j)}\cr}$$
On the right, the first summation is
\nobreak
$$\sum_{a\in G}\mu[a,g]\sigma_2(a)\mu[(i,j);a]=
  \delta\ce(i)\delta_{i,\omega}\delta_{i,\pi(j)}$$
The second summation is
\nobreak
$$\sum_{r,s,t}\mu[(r,s);g]\sigma_4[t,r]\mu[(i,j);(t,s)]$$
Note this will be nonzero iff there exist indices $r,s,$ and $t$ satisfying
$$\eqalign{r&=s\ci=\pi(s)\inv\cr
t&=\pi(r\inv)\cr
i\circ j&=t\cr
t*s&=i\cr}$$
We can solve for the indices uniquely, in terms of $i$ if they exist:
$$\eqalign{ &\pi\inv(s)*s=i\cr
\hbox{so}\quad&\pi(s)\inv*\pi(s)*\pi\inv(s)*s=i\cr
\hbox{and}\quad &s=\pi\inv(\omega*i\inv)\cr}$$
Note that if $i=\omega$ there could be no $s$ filling the requirement
and the sum would vanish.

Immediately, $t=\pi(\omega*i\inv)$ follows.  Now from the remaining equation $i\circ j = t$ we discover
the condition necessary to get a nonzero sum:

$$\eqalign{j=t\circ i\ci &= \pi\bigl(\pi\inv(t)*\pi\inv(\pi\inv(i\inv))\bigr)\cr
&=\pi\bigl(\omega*i\inv*\pi(i\inv)\bigr)\cr
&=\pi\bigl(\omega*i\inv*\pi(i\inv)*\pi\inv(i\inv)*\pi\inv(i\inv)\inv\bigr)\cr
&=\pi(\pi(i))=\pi\inv(i)\cr}$$
So a nonzero summation happens precisely when $i=\pi(j)$ and the reduction is
$$\eqalign{\sum_{r,s,t}&\mu[(r,s);g]\sigma_4[t,r]\mu[(i,j);(t,s)]\cr
&\qquad=\ce(i)\chi_\omega(g)\inv\chi_i(g)\inv\psi(\omega*i)
  N\left(\pi(\omega i\inv),\pi\inv(\omega i\inv)\right)
  (1-\delta_{i,\omega})\delta_{i,\pi(j)}\cr}$$
The upshot is that for $|G|$ odd, $\omega=\epsilon$, this piece of the
pentagon is equivalent to the requirement
$$\sigma_3(\epsilon)\psi(\pi(i)\inv)=\psi(i)
  N\left(\pi(i\inv),\pi\inv(i\inv)\right),\quad\hbox{for all $i$}$$
For $|G|$ even, $\omega\ne\epsilon$, this piece of the pentagon is
equivalent to the requirements
$$\eqalign{\sigma_3(\epsilon)\psi(\pi\inv(\omega))&=1\cr
  \sigma_3(\epsilon)\psi(\pi(i)\inv)&=\psi(\omega*i)
  N\left(\pi(\omega*i\inv),\pi\inv
  (\omega*i\inv)\right),\quad\hbox{for all\ }i\ne\omega\cr}$$

\sh{8.10\qua Reducing the $mmm/m$ hexagon, part III }

This is very similar to Part II, so we omit some details. On the left,
$$\eqalign{&\sum_t \sigma_4[t,s]\mu[g;(r,t)\sigma_3(g)\cr
&\qquad=\sigma_3(\epsilon)\chi_\omega(g)\psi(s)
\left(|G|\xi\left(\pi(s)\right)\chi_{\pi(s)}(g)\inv
  \ce\left(\pi(s)\right)\right)\inv \delta_{s,\pi(r)}\cr}$$
On the right, the first summation is
$$\sum_{a\in G}\mu[a;(r,s)]\sigma_2(a)\mu[g,a]=
  \delta\left(|G|\xi(\omega)\ce(\omega)\right)\inv
  \delta_{r,s\inv}\delta_{\omega,\pi(s)}$$
The second summation is 
$$\eqalign{&\sum_{j,p,q}\mu[(p,q);(r,s)]\sigma_4[j,p]\mu[g;(j,q)]\cr
&\qquad=r_g\left(\pi\left(r*s\right)\inv\right)
  \psi(r*s)N(r,s)\delta_{s,\pi(r)}(1-\delta_{r,s\inv})\cr
&\qquad=\left(|G|\xi\left(\omega*\pi(s)\right)
\chi_{\omega*\pi(s)}(g)\inv
\ce\left(\omega*\pi(s)\right)\right)\inv
\cr&\qquad\qquad\qquad\cdot
\psi\left(\omega*\pi(s)\inv\right)
N\left(\pi\inv(s),s\right)
\delta_{s,\pi(r)}(1-\delta_{r,s\inv})\cr}$$
So for $|G|$ odd, $\omega=\epsilon$, the requirement of the pentagon is 
$$\sigma_3(\epsilon)\psi(s)=\psi\left(\pi(s)\inv)\right)
  N\left(\pi\inv(s),s\right),\quad\hbox{for all $s$}$$
Note that this is equivalent to the requirement from Part II of the pentagon, under the change of variable $s\to\pi(i)\inv$.

For $|G|$ even, $\omega\ne\epsilon$, this piece of the pentagon is
equivalent to the requirements
$$\sigma_3(\epsilon)\psi\left(\pi\inv(\omega)\right)=1$$
and
$$\eqalign{&\sigma_3(\epsilon)\psi(s)\xi\left(\omega*\pi(s)\right)
  \ce\left(\omega*\pi(s)\right)\cr
&\qquad=\xi\left(\pi(s)\right)\ce\left(\pi(s)\right)\psi
  \left(\omega*\pi(s)\inv\right)N\left(\pi\inv(s),s\right)\cr}$$
for all $s\ne\pi\inv(\omega)$.  Note that while the first equation
duplicates the result from Part II, the second is a bit more complicated,
as it involves $\xi$'s.

\sh{8.11\qua Reducing the $mmm/m$ hexagon, part IV}

On the left,
$$\eqalign{&\smash{\sum_{u,v}}\sigma_4[u,s]\mu[(i,v);(r,u)]\sigma_4[j,v]\cr
&\qquad=\sigma_4[\pi(s\inv),s]
  \mu[(i,\pi(j\inv));(r,\pi(s\inv))]\sigma_4[j,\pi(j\inv)]\cr
&\qquad=\psi(s)\psi\left(\pi(j\inv)\right)N\left(r,\pi(s\inv)\right)
\delta_{i,r*\pi(s\inv)}\delta_{r,i\circ\pi(j\inv)}\cr}$$
On the right, the first summation is
$$\eqalign{&\sum_{a\in G}\mu[a;(r,s)]\sigma_2(a)\mu[(i,j);a]\cr
&\quad=\left(|G|\xi\left(\pi\inv(r)\right)
  \ce\left(\pi\inv(r)\right)\right)\inv\ce(i)
  \sum_a\chi_\omega(a)\chi_i(a)\inv
  \chi_{\pi\inv(r)}(a)\inv\delta_{r,s\inv}\delta_{i,j\ci}\cr
&\quad=\left(\xi\left(\pi\inv(r)\right)\ce\left(\pi\inv(r)\right)\right)\inv
\ce(i)\delta_{r,s\inv}\delta_{i,j\ci}\delta_{\omega,i*\pi\inv(r)}\cr}$$
Recall the second summation is
\nobreak
$$\sum_{p,q,t}\mu[(p,q);(r,s)]\sigma_4[t,p]\mu[(i,j);(t,q)]$$
As usual, there could be at most one nonzero term in the sum and the
game is to find the conditions on $i,j,r,s$ which would permit this.

A nonzero sum happens iff there are $p,q,t$ with
\medskip
\centerline{\vbox{\halign{ $#$&$#$&$#$\cr
r=p\circ q&&p=r*s\cr
t=i\circ j&&i=t*q\cr
&t=\pi(p\inv)&\cr}}}
\medskip
In terms of $r$ and $s$, the magic indices must be
$$\eqalign{p&=r*s\cr
q&=\pi(s)*\pi(r*s)\inv\cr
t&=\pi\left(\left(r*s\right)\inv\right)\cr}$$
but in terms of $i$ and $j$ we must have
$$\eqalign{p&=\pi(i\inv)*\pi(j\inv)\cr
q&=\pi\inv\left(j*\pi(i)\right)\cr
t&=i\circ j = \pi\left(\pi\inv(i)*\pi\inv(j)\right)\cr}$$
So we can determine $r$ and $s$ in terms of $i$ and $j$ if there is to
be a nonzero sum.

First, we can solve for $s$:
$$\eqalign{\pi(s)&=\pi\inv(j*\pi(i))*\pi(r*s)\cr
&=\left(\pi(j\inv*\pi(i)\inv)\right)\inv*
  \pi\left(\pi(i\inv)*\pi(j\inv)\right)\cr
&=\left(i\inv*\pi\left(\pi(j\inv)\inv*
  \pi\inv(i)\right)\right)\inv*\pi\left(\pi(i\inv)*\pi(j\inv)\right)\cr
&=i*\pi\inv\left(\pi(j\inv)*\pi(i\inv)\right)*
  \pi\left(\pi(i\inv)*\pi(j\inv)\right)\cr
&=i*\omega*\pi\inv(j)*\pi\inv(i)\cr
&=\omega*\omega*\pi(i)\inv*\pi\inv(j)\cr
&=\pi(i)\inv*\pi\inv(j)\cr}$$
and the conclusion is
$$s = \pi\inv\left(\pi(i)\inv*\pi\inv(j)\right)=
  \pi\left(\pi(i)*\pi(j\inv)\right)\inv$$
Now we can solve for $r$ as well:
$$\eqalign{r&=\pi(i\inv)*\pi(j\inv)*\pi\left(\pi(i)*\pi(j\inv)\right)\cr
&=\pi(i\inv)*\pi(j\inv)*\pi\inv(i)*\pi\left(\pi(j)*i\inv\right)\cr
&=\pi(j\inv)*\pi\inv(j)*\pi\left(\pi\inv(i)*j\inv\right)\cr
&=\pi\left(\pi\inv(i)*j\inv\right)\cr}$$
and so the complete reduction of the second summation is:
\nobreak
$$\eqalign{\smash{\sum_{p,q,t}}&\mu[(p,q);(r,s)]\sigma_4[t,p]\mu[(i,j);(t,q)]\cr
&\quad=N\left(r,s\right)\psi(r*s)N
  \left(\pi\left((r*s)\inv\right),\pi(s)*\pi(r*s)\inv\right)\cr
&\quad\phantom{=}\quad\delta_{r,\pi(\pi\inv(i)*j\inv)}
  \delta_{s,\pi\inv(\pi(i)\inv)*\pi\inv(j))}
  (1-\delta_{r,s\inv})(1-\delta_{i,j\ci})\cr}$$
Now we can set about extracting practical equations.

Assume the lefthand side is nonzero, ie, $i=r*\pi(s\inv)$ and
$r=i\circ\pi(j\inv)$.

{\bf Case I}\qua  $r=s\inv$ and $i=j\ci$\qua
In this case the second sum on the right vanishes but since we have
$i=r*\pi(r)=\omega*\pi(r\inv)$, the first sum on the right does not
vanish, so we get the requirement
$$\psi(s)\psi(\pi(j\inv))N(r,\pi(s\inv))=
  \left(\xi(\pi\inv(r))\ce(\pi\inv(r))\right)\inv\ce(i)$$
This can be rewritten all in terms of $r$, but it requires some work to
express $\pi(j\inv)$ in terms of $r$:
$$\eqalign{\pi(j\inv)&=r\circ\left(r*\pi(r)\right)\ci\cr
&=\pi\left(\pi\inv(r)*\pi(r\inv*\pi(r)\inv)\right)\cr
&=\left(\pi(r*\pi(r))\right)\inv*\pi\inv(r)\cr
&=\pi\inv(r*\pi(r))\cr}$$
So our final version of the hexagon in this case is
$$\psi(r\inv)\psi(\pi\inv(r*\pi(r)))N(r,\pi(r))
=\left(\xi(\pi\inv(r)\ce(\pi\inv(r))\right)\inv\ce(r*\pi(r)),$$
for all $r$.

{\bf Case II}\qua $r\ne s\inv$ and $i\ne j\ci$\qua
In this case the first sum on the right vanishes but it is easy to verify
that the second does not.  As in the previous case we can
solve for $\pi(j\inv)$ in terms of $r$ and $s$, namely,
$\pi(j\inv)=\pi\inv(s\inv*\pi(r))$.  So our final version of the hexagon
in this case is
$$\eqalign{&\psi(s)\psi\left(\pi\inv(s\inv*\pi(r))\right)N(r,\pi(s\inv))\cr
&\qquad=N(r,s)\psi(r*s)N
  \left(\pi\left((r*s)\inv\right),\pi(s)*\pi(r*s)\inv\right)\cr}$$
for all $r\ne s\inv$.

{\bf Case III}\qua $r\ne s\inv$ xor $i\ne j\ci$\qua
If the lefthand side is nonzero, it is easy to check that this cannot occur.

That takes care of all possibilities when the lefthand side is nonzero.
But conversely any conditions under which the righthand side is nonzero
give a nonzero lefthand side so we have covered all the nontrivial
possibilities, and completed reducing all the hexagon equations.  As with
the results of the largest pentagon in section 5, the expressions are
awkward but do have the virtue of admitting an obvious trivial solution.

We summarize the work of this section so far in the following:

\proclaim {8.12\qua Proposition}
Suppose near-group fusion rule $(G,k)$, $|G|=k+1$, has a monoidal
structure described by $\xi,\ce$, and $N$.  If a braiding is possible,
then that braiding has the following structure:
\par{\rm(1)}\qua $\sigma_0\equiv 1$
\par{\rm(2)}\qua $\sigma_1(g)=\sigma_2(g)=\chi_\omega(g)$ for all $g\in G$.
\par{\rm(3)}\qua $\sigma_3(g)=\sigma_3(\epsilon)\chi_\omega(g)$ for all $g$ in $G$
\par{\rm(4)}\qua $\sigma_4(g)=\bigl(\psi(j)\delta_{i,\pi(j\inv)}\bigr)_{i,j}$
where $\sigma_3(\epsilon)$ and $\psi(1),\ldots,\psi(k)$ are invertible
constants.
\endproclaim

\prf
Summary of results up to this point.
\endprf

If we want to construct a braiding, therefore, 8.12 gives a recipe for the commuting maps
in terms of constants $\sigma_3(\epsilon)$ and $\psi(1),\ldots,\psi(k)$.  The following summarizes
the constraints on those constants.

\proclaim {8.13\qua Proposition} The hexagon axiom is satisfied provided the constants $\sigma_3(\epsilon)$ and $\psi(1),\ldots,\psi(k)$
satisfy the following constraints:\endproclaim

Independent of $|G|$, for all indices $r$,
$$\eqalignno{
  \psi(r)\psi(\pi(r)\inv)\xi(\pi(r\inv))&=\sigma_3(\epsilon)
  \xi(r)\xi(\pi(r))&(29)\cr}$$
\indent for $r\ne\pi(r)\inv$,
$$\eqalignno{\psi(r\inv)\psi(\pi\inv(r*\pi(r)))N
  (r,\pi(r))&\hfill&\cr
&\hskip-1in=\left(\xi(\pi\inv(r)\ce(\pi\inv(r))\right)\inv
  \ce(r*\pi(r))&(30)\cr}$$
\indent and for $s\ne\pi(r)$, $s\ne r\inv$,
$$\eqalignno{\psi(s)\psi\left(\pi\inv(s\inv*\pi(r))\right)N
  (r,\pi(s\inv))&\hfill&\cr
&\hskip-1.125in=N(r,s)\psi(r*s)N
  \left(\pi\left((r*s)\inv\right),\pi(s)*\pi(r*s)\inv\right)&(31)\cr}$$
For $|G|$ odd, for all indices $r$,
$$\eqalignno{
\sigma_3(\epsilon)^2&=\delta&(32)\cr
\sigma_3(\epsilon)\psi(\pi(r)\inv)&=
  \psi(r)N\left(\pi(r\inv),\pi\inv(r\inv)\right)&(33)\cr}$$
For $|G|$ even,
$$\eqalignno{
\sigma_3(\epsilon)^2&=\psi(\omega)\xi(\omega)\inv&(34)\cr
\sigma_3(\epsilon)\psi(\pi\inv(\omega))&=1&(35)\cr
\sigma_3(\epsilon)\psi(\pi(r)\inv)&=
  \psi(\omega*r)N\left(\pi(\omega*r\inv),\pi\inv(\omega*r\inv)\right),\quad
  \hbox{for\ }r\ne\omega&(36)\cr}$$
$$\eqalignno{
\sigma_3(\epsilon)\psi(\pi\inv(r))\xi\left(\omega*r\right)\ce
  \left(\omega*r\right)
&\hfill&\cr
&\hskip-1.125in=\xi\left(r\right)\ce\left(r\right)\psi\left(\omega*r\inv\right)
  N\left(\pi(r),\pi\inv(r)\right),\quad\hbox{for\ }r\ne\omega&(37)\cr}$$

\prf
Equations (29)--(31) come from 8.7 and 8.11.  Equations (32) and (33)
come from 8.8 and 8.9.  Equations (34)--(37) come from 8.8, 8.9, and 8.10,
with minor change of indices to make the present list appear consistent.
\endprf

\proclaim {8.14\qua Proposition}
The standard monoidal solution constructed in 5.18 admits a symmetric
commuting structure.
\endproclaim

\prf
This is now a triviality.  The standard monoidal solution
corresponds to taking $\xi,\ce$, and $N$ all identically
equal to 1.  Now we can satisfy the constraints listed in
8.3 by setting $\sigma_3(\epsilon)=\psi(1)=\cdots=\psi(k)=1$
as well, obtaining a solution to the hexagon equations.
It's easy to check that this is a symmetric structure:
$\sigma_0(g,h)\sigma_0(h,g)=1$ since $\sigma_0\equiv 1$.  Also
$\sigma_1(g)\sigma_2(g)=\sigma_2(g)\sigma_1(g)=\chi_\omega(g)^2=1$.
Similarly, $\sigma_3(g)^2=\sigma_3(\epsilon)\chi_\omega(g)^2=1$, and
$\sigma_4=\left(\delta_{i,\pi(j\inv)}\right)_{i,j}$ is self-inverse
because the map $j\mapsto\pi(j\inv)$ is an involution.  Since the
structure is symmetric we do not need to consider any inverse hexagons.
\endprf

\section{Examples of braidings}

\sh {9.1\qua Example 1\qua $k=1, G={\Bbb Z}/2{\Bbb Z}$ }

From section 3, there are three monoidal structures, corresponding to a
choice of $\xi$ a cube root of unity.  In the language we have developed
the monoidal structure is described by $\ce(1)=1$, $xi(1)=\xi$, and
there is no $N$ due to a shortage of indices.  Braidings, if possible,
are described by invertible constants $\sigma_3(\epsilon)$ and $\psi(1)$.
The possibilities are given in the following:

\proclaim {9.2\qua Theorem {\rm (Braidings, $k=1, G={\Bbb Z}/2{\Bbb Z}$)}}
The $\xi=1$ structure admits three distinct braidings corresponding to
a choice of $\psi=\psi(1)$ with $\psi^3=1$.  The $\xi\ne1$ structures
do not admit braidings.
\endproclaim

\prf
Referring to proposition 7.13 we see that the only requirements of the hexagon
axiom are:
$$\eqalign{\psi^2\xi&=\sigma_3(\epsilon)\xi^2\cr
\sigma_3(\epsilon)^2&=\psi\xi\inv\cr
\sigma_3(\epsilon)\psi&=1\cr}$$
(these come from equations (29), (34), and (35), respectively;
none of the remaining equations from 8.13 apply in this case).
It is easy that these equations are satisfied iff $\psi^3=\xi$ and
$\sigma_3(\epsilon)=\psi\inv$.  But inverse hexagons will be satisfied iff
$\psi^{-3}=\xi$ as well; we conclude that $\xi=1$ is the only monoidal
structure which admits any commutative structure satisfying both the
hexagon and inverse hexagon axioms.  Of course the $\psi=1$ structure
is symmetric; the other two clearly not.

For reference, then, the complete data for the braidings in this case is:
$$\eqalign{\sigma_0&\equiv (1)\cr
\sigma_1(a)=\sigma_2(a)&=\pmatrix {\chi(a)} \cr
\sigma_3(a)&=\pmatrix {\psi\inv\chi(a)}\cr
\sigma_4&=\pmatrix{\psi}\cr}$$
where $\chi$ is the nontrivial character of $G$ and $a\in G$.
\endprf

\sh {9.3\qua Example 2\qua $k=2, G={\Bbb Z}/3{\Bbb Z} = \{\epsilon,g,g^2\}$ }
From section 3 there are two distinct monoidal structures corresponding to
a choice of $\xi=\pm1$.  In the language we have developed, the monoidal
structure is described by
$$\eqalign{\xi(1)=&\xi(2)=\xi\cr
  \ce(1)=1,\quad&\ce(2)=\xi\cr
  N(1,1)=1,\quad&N(2,2)=\xi\cr}$$
Braidings, if possible, are described by $\sigma_3(\epsilon)$, $\psi(1)$,
and $\psi(2)$. The possibilities are given in the following:

\proclaim {9.4\qua Theorem {\rm (Braidings, $k=2, G={\Bbb Z}/3{\Bbb Z}$)}}
The $\xi=1$ structure admits four distinct braidings corresponding to a
choice of $\psi(1),\psi(2)$ with $\psi(1)^2=\psi(2)^2=1$.  The $\xi=-1$
structure does not admit any braiding.
\endproclaim

\prf
Referring to proposition 8.13, and using the monoidal data, we see the
requirements of the hexagon axiom are as follows.  From equation (29),
$$\psi(1)\psi(2)\xi=\sigma_3(\epsilon)$$
From equation (30),
$$\eqalign{\psi(2)^2&=1\cr
  \psi(1)^2\xi&=1\cr}$$
There are no indices $r,s$ to which (31) applies.  From (32),
$$\sigma_3(\epsilon)^2=\xi$$
and from (33),
$$\eqalign{\sigma_3(\epsilon)\psi(2)&=\psi(1)\xi\cr
\sigma_3(\epsilon)\psi(1)&=\psi(2)\cr}$$
It's clear that these are satisfied iff we take $\psi(1)^2=\xi$,
$\psi(2)^2=1$, and $\sigma_3(\epsilon)=\psi(1)\psi(2)\xi$.  For inverse
$mmm/m$ hexagon to be satisfied, though, $\sigma_3(\epsilon)$ has
to be self-inverse, and so $\xi=+1$ is the only possibility for the
monoidal structure. The cases $\psi(1)=\psi(2)=\pm1$ are symmetric;
the $\psi(1)\ne\psi(2)$ cases are not.

For reference, the complete data for the braidings in this case is:
$$\eqalign{\sigma_0&\equiv(1)\cr
  \sigma_1\equiv\sigma_2&\equiv(1)\cr
  \sigma_3(a)&=\pmatrix {\psi(1)\psi(2) \cr}\cr
  \sigma_4&=\pmatrix {0&\psi(2) \cr \psi(1)&0 \cr}\cr}$$
where $a\in G$.
\endprf

\sh {9.5\qua Example 3\qua
     $k=3, G={\Bbb Z}/4{\Bbb Z} = \{\epsilon,g,g^2,g^3\}$ }

From section 3, this fusion rule has a unique monoidal structure,
which is described by taking $\xi,\ce$, and $N$ all identically 1.
As it turns out, the commutative structure is unique as well.

\proclaim {9.6\qua Theorem {\rm(Braidings, $k=3, G={\Bbb Z}/4{\Bbb Z}$)}}
The unique monoidal structure for this fusion rule admits a unique
braiding.
\endproclaim

\prf
Referring to proposition 8.13, and using the monoidal data, we get the
following constraints.  From (29) with $r=1,2,3$, respectively:
\eject
$$\eqalign{\psi(1)\psi(2)=\sigma_3(\epsilon)\quad&
  \quad\psi(2)\psi(1)=\sigma_3(\epsilon)\cr
  \psi(3)^2&=\sigma_3(\epsilon)\cr}$$
From (30) with $r=1,2$:
$$\psi(3)\psi(2)=1\quad\quad\psi(2)\psi(3)=1$$
From (31) with $(r,s)=(2,1),(3,2)$:
$$\psi(1)^2=\psi(3)\quad\quad\psi(2)^2=\psi(1)$$
From (34) and (35),
$$\sigma_3(\epsilon)^2=\psi(2)\quad\quad\sigma_3(\epsilon)\psi(1)=1$$
And from (36),
$$\sigma_3(\epsilon)\psi(2)=\psi(3)\quad
  \quad\sigma_3(\epsilon)\psi(3)=\psi(1)$$
Since the monoidal data is trivial, equation (37) is equivalent to
(36).  It is easy to check that the solutions to these are of the
following form: $\psi(2)$ any fifth root of unity; $\psi(1)=\psi(2)^2$;
$\psi(3)=\psi(2)\inv$; $\sigma_3(\epsilon)=\psi(2)^{-2}$.  In other
words, the commutativity data is determined by the choice of $\psi(2)$.
However, inverse hexagons (consider $mmm/\epsilon$) force the additional
requirement $\psi(1)=\psi(2)=\psi(3)$, so the only braiding comes
from $\psi(2)=1$ (our standard symmetric solution).  For reference,
the complete data:
$$\eqalign{\sigma_0&\equiv(1)\cr
  \sigma_1(a)=\sigma_2(a)&=\pmatrix {\chi_\omega(a)} \cr
  \sigma_3(a)&=\pmatrix {\chi_\omega(a)}\cr
  \sigma_4&=\pmatrix {1&0&0 \cr 0&0&1 \cr 0&1&0 \cr}\cr}$$
where $\chi_\omega$ is the nontrivial order 2 character and $a\in G$.
\endprf

\section {Additional structures II:  twist} %

This is just a brief comment on the possibility of adding twist morphisms
to a braided near-group category in the setting of Theorem 1.2 (useful
for invariants of framed links).  We need endomorphisms $\theta_s$
for each simple $s$, which balance the commutative structure in the
sense that $\theta_{r\tens s}=\theta_r\theta_s\sigma(r,s)\sigma(s,r)$
where $\sigma(r,s)$ is used generically to denote commuting $r$ past $s$.
After the simplifications of the braided structure carried out in section
8, it is easy to see that this axiom reduces to the following:
\eject

$$\eqalign{\theta_g&=1\quad\hbox{for all\ }g\in G\cr
\sigma_3(\epsilon)^2\theta_m^2&=1\cr
\psi(j)\psi\left(\pi(j\inv)\right)&=\theta_m\quad\hbox{for all\ }j\cr}$$
This is simple to work out for the small example categories we have
studied in sections 3 and 9, but interesting since it gives a few
examples of non-symmetric braidings that do or do not admit braidings.
We say a braiding is balanced if there is a $\theta_m$ satisfying the
twist equations above.

\proclaim {10.1\qua Proposition {\rm (Twist, $k=1, G={\Bbb Z}/2{\Bbb Z}$)}}
In the setting of theorem 9.2, the $\psi=1$ structure is balanced.
The $\psi\ne1$ structures are not.
\endproclaim

\prf
The symmetric $\psi=1$ case has trivial twist morphisms and the braiding data
from 8.2 together with the twist equations above implies $\psi=1$ immediately.
\endprf

\proclaim {10.2\qua Proposition {\rm (Twist, $k=2, G={\Bbb Z}/3{\Bbb Z}$)}}
In the setting of theorem 9.4, all of the braidings are
balanced.
\endproclaim

\prf
The symmetric $\psi(1)=\psi(2)$ structures take trivial twist morphisms;
taking $\theta_m=-1$ satisfies the nonsymmetric $\psi(1)=-\psi(2)$
structures.
\endprf
\bigskip
{\large\bf References}\bigskip

{\small\leftskip 30pt

\item{[FK]} {\bf J Frolich}, {\bf T Kerler},
  {\it Quantum Groups, Quantum Categories, and quantum field theory},
  Springer--Verlag (1993)
\item{[KR]} {\bf R M Kashaev}, {\bf N Reshitikhin},
  {\it Symmetrically Factorizable Groups and Set--Theoretical Solutions
  of the Pentagon Equation},
  {\tt math.qa/0111171}
\item{[Mac]} {\bf S MacLane},
  {\it Categories for the Working Mathematician},\nl
  Springer--Verlag (1971)
\item{[O]} {\bf V Ostrik},
  {\it Fusion categories of rank 2},
  {\tt math.qa/0203255}
\item{[Q]} {\bf F Quinn},
  {\it Lectures on Axiomatic Topological Quantum Field Theory},
  Geometry and Quantum Field Theory Volume 1,
  IAS/Park City Mathematical Series,
  American Mathematical Society (1995)
\item{[TY]} {\bf D Tambara}, {\bf S Yamagami},
  {\it Tensor Categories with Fusion Rules of Self-Duality for Finite
  Abelian Groups},
  Journal of Algebra 209 (1998), 692--707
\item{[S]} {\bf J Siehler},
  {\it Braided Near-group Categories},
  {\tt math.qa/0011037}
\enditems
}
\medskip
\Addressesr

\bye